%% file: main.tex
\documentclass[11pt]{article}
\usepackage{fullpage}

\input{macros}

\input{tempmacros}

\begin{document}

\title{Metastable Mixing of Markov Chains: Efficiently Sampling \\ Low Temperature Exponential Random Graphs%\thanks{Grants or other notes
}
% \subtitle{}

% \titlerunning{Sampling Low-Temperature Exponential Random Graphs}        % if too long for running head

% \author{Guy Bresler \and
% 			Dheeraj Nagaraj       \and
%         Eshaan Nichani %etc.
% }

\author{Guy Bresler \\ \texttt{MIT}\\\texttt{guy@mit.edu} \and 
			Dheeraj Nagaraj  \\ \texttt{Google Research} \\ \texttt{dheerajnagaraj@google.com}    \and
        Eshaan Nichani \\ \texttt{Princeton University} \\ \texttt{eshnich@princeton.edu}%etc.
}


% \institute{G. Bresler \at
%               %first address \\
%               %Tel.: +123-45-678910\\
%               %Fax: +123-45-678910\\
%               \email{guy@mit.edu}           %  \\
% %             \emph{Present address:} of F. Author  %  if needed
%           D. Nagaraj \at
%              \email{dheerajnagaraj@google.com}
%           E. Nichani \at
%              \email{eshnich@princeton.edu}
% }


\maketitle

\begin{abstract}
% The notion of mixing for finite state Markov chains is usually defined for worst-case initial distribution, 
% and this can be overly restrictive for many well studied examples. 
% In this paper we consider the problem of 
% sampling from the low-temperature exponential random graph model (ERGM). The usual approach is via Markov chain Monte Carlo, but strong lower bounds have been established for the ERGM showing that any local Markov chain suffers from an exponentially large mixing time due to metastable states.
% We instead consider \emph{\stable~mixing}, a notion of approximate mixing within a collection of metastable states.
% % convergence of a Markov chain to the stationary distribution up to a total variation distance of $\delta$ with a given initial distribution. 
% In the case of the ERGM, we show that Glauber dynamics with the right 
% % (mixture of) 
% $G(n,p)$ initialization has a \stable~mixing time of $O(n^2\log n)$ to within total variation distance $\exp(-\Omega(n))$. 
% 
% 
In this paper we consider the problem of 
sampling from the low-temperature exponential random graph model (ERGM). The usual approach is via Markov chain Monte Carlo, but 
% strong lower bounds have been established 
Bhamidi et al. showed that any local Markov chain suffers from an exponentially large mixing time due to metastable states.
We instead consider \emph{\stable~mixing}, a notion of approximate mixing relative to the stationary distribution, for which it turns out to suffice to mix only within a collection of metastable states.
% convergence of a Markov chain to the stationary distribution up to a total variation distance of $\delta$ with a given initial distribution. 
We show that the Glauber dynamics for the ERGM at any temperature -- except at a lower-dimensional critical set of parameters -- when initialized at
% (mixture of) 
$G(n,p)$ for the right choice of $p$ has a \stable~mixing time of $O(n^2\log n)$ to within total variation distance $\exp(-\Omega(n))$. 

% When the stationary distribution is the low temperature exponential random graph model (ERGM), strong lower bounds have been established which show the Markov chain suffers from an exponentially large mixing time whenever the Markov chain is local. In this work, we consider \emph{\stable~mixing}, i.e., convergence of a Markov chain to the stationary distribution up to a total variation distance of $\delta$ with a given initial distribution. In the case of ERGM, we show that Glauber dynamics with the right (mixture of) $G(n,p)$ initialization has an \stable~mixing time of $O(n^2\log n)$ when the considered TV distance $\delta \geq \exp(-\Omega(n))$. 
%We also show a multi timescale behavior (called \textbf{\slack}) associated with \stable~mixing which is in a sense the opposite of the well studied cut-off phenomenon. In the ERGM case, this means that whenever $\delta \leq \exp(-\Omega(n^2))$, the \stable~mixing time for the same initialization is $\exp(\Omega(n))$.

%Our main technical result combines the coarse grained large deviations results for ERGM in the graphon space with the cavity method to obtain fine grained uniform concentration results for the sub-graph counts of ERGM. This result, in particular, establishes fast \stable~mixing of Glauber dynamics with the right initialization.
%\keywords{Random Graphs \and Exponential Random Graphs \and Markov Chain Monte Carlo}
\end{abstract}

\section{Introduction}
\label{intro}

% In this work, we consider the problem of sampling from the exponential random graphs model (ERGM). 

Given a vector of real-valued parameters 
$\beta := (\beta_0,\beta_1,\dots,\beta_K) \in \mathbb{R}\times (\mathbb{R}^{+})^{K}$
% $\beta := (\beta_0,\beta_1,\dots,\beta_K)$
, the exponential random graph model $\ergm(n,\beta)$ is defined to be the probability measure over all simple graphs with $n$ vertices 
$$\mu_\beta(X) = \frac1{Z_\beta}\exp\bigg(  \sum_{i=0}^{K}n^2\beta_i N_i(X)\bigg)\,.$$
The $N_i$ are subgraph counts corresponding to finite graphs $G_0,G_1,\dots, G_K$ (such as edges, triangles, 4-cycles, 5-cycles, 2-stars, etc.) and $Z(\beta)$ denotes the normalizing constant.
This is an exponential family where the sufficient statistics are the subgraph counts. The model significantly generalizes the Erd\H{o}s-R\'enyi random graph and is used to model a variety of complex networks like social networks and biological networks \cite{frank1986markov,holland1981exponential, fienberg2010introduction,fienberg2010introduction2,wasserman1994social}. Early analysis was carried out by statistical physicists \cite{park2004solution,park2005solution,burda2004network}, and probabilists and statisticians have further studied various questions about these models including sampling, estimation, large deviations theory, concentration of measure, and phase transitions \cite{bhamidi2011mixing,chatterjee2007stein,radin2013phase,eldan2017exponential,ganguly2019sub,mukherjee2013statistics,yin2017asymptotics,reinert2019approximating}. 

The basic problem we consider in this work is that of producing a sample from the $\ergm$ probability distribution in polynomial time. A popular approach to sampling is to use the Glauber dynamics, a simple reversible Markov chain with the desired stationary distribution, and to run it for sufficiently long that it is close to stationarity.

\begin{definition}[Glauber Dynamics]
\label{def:glauber}
Given any probability distribution $\pi$ over $\Omega = \{0,1\}^{{{n}\choose{2}}}$ with $\pi(x)>0$ for all $x\in \Omega$, the Glauber dynamics with respect to $\pi$ is the discrete time Markov chain over $\Omega$ with single step transition from $X$ to $X^{\prime}$ as follows:
\begin{enumerate}
    \item Pick a coordinate $E\in\left[{{n}\choose{2}}\right]$ uniformly at random. 
    \item Form $X'$ by keeping all edges except for $E$ the same as in $X$ and sample $X'_E\sim \pi(\cdot|X_{\sim E})$ conditional on the other edges.
\end{enumerate}

 The Glauber dynamics with respect to $\pi$ is reversible for $\pi$ (\cite{levin2017markov}), so in particular $\pi$ is stationary.
\end{definition}

The mixing time of the Glauber dynamics, i.e., the time until the distribution is within total variation $1/4$ of stationarity, determines whether the approach is feasible and as shown by Bhamidi et al.~\cite{bhamidi2011mixing} turns out to be essentially characterized by
the function 
$L_{\beta} : [0,1] \to \mathbb{R}$ defined as 
\begin{equation}
\label{e:Lbeta}
    L_{\beta}(p) = \sum_{i=0}^{K}\beta_i p^{|E_i|} - I(p)
\end{equation}
where $I(p) := \frac{1}{2}p\log p+\frac{1}{2}(1-p)\log (1-p)$. (Bhamidi et al.~\cite{bhamidi2011mixing} actually formulated their results in an equivalent way in terms of the function $\phi_\beta$ defined in Section~\ref{sec:Glauber}, while $L_\beta$ was studied by \cite{chatterjee2013estimating}.)  

\begin{theorem}[\cite{bhamidi2011mixing}]
Consider the $\ergm(n,\beta)$ distribution. There are three regimes for $\beta$:
\begin{enumerate}
    \item High temperature: If $L_{\beta}$ has a unique local maximum with a non-vanishing second derivative, then the Glauber dynamics Markov chain for $\ergm(n,\beta)$ mixes in time $O(n^2 \log n)$.
    \item  Low temperature: If $L_\beta$ has multiple local maxima with non-vanishing second derivatives, then any local Markov chain with stationary distribution $\ergm(n,\beta)$ must suffer a mixing time of $\exp(\Omega(n))$. 
    \item Critical temperature: If any local maxima of $L_{\beta}$ has a vanishing second derivative.
\end{enumerate}
\end{theorem}
In this work, we use a loose and intuitive notion of metastability since this is not important in order to state our technical results:  consider subsets of the state space $A,B \subseteq \Omega$ such that $A \subseteq B$. We will call $B$ to be metastable with respect to a given markov dynamics if the markov dynamics initialized inside the set $A$ takes a long time to exit the set $B$. In this work, we interpret `long time' as being exponential in $n$. Slow mixing in the low-temperature phase is due to the existence of multiple, disconnected metastable states from which it takes the Glauber dynamics exponential time to leave. The question is therefore: can one efficiently produce a sample from the $\ergm$ in the low temperature regime?

An important insight into the structure of the $\ergm$ distribution was developed by Chatterjee and Diaconis \cite{chatterjee2013estimating}, and this will constitute a useful step towards our goal. 
They showed that the $\ergm$ distribution is close to a finite mixture of constant graphons with respect to the cut-metric:

\begin{theorem}[Theorem 4.2 of \cite{chatterjee2013estimating}]
\label{thm:large_deviations} Denote by $M_{\beta}$ the set of global maxima of $L_\beta$. The $\ergm(n,\beta)$ distribution converges, in probability with respect to the cut-metric, to a mixture of $G(n,p^{*})$ for $p^{*} \in M_{\beta}$.
Formally, let $X_{n} \sim \ergm(n,\beta)$ and $\tilde{X}_n$ be its corresponding graphon and $\tilde{M}_{\beta}$ be the set of all constant graphons with value $p^{*}$ for some $p^{*} \in M_{\beta}$. For every fixed $\eta > 0$, there are constants $C(\eta),c(\eta) > 0$ such that
$$\mathbb{P}(\delta_\square(\tilde{X}_n,\tilde{M}_{\beta})> \eta) \leq C(\eta)\exp(-c(\eta)n^2)\,.$$
\end{theorem}

Graphons and the cut metric are reviewed in Section~\ref{s:graphon}.

Given the approximation results of Chatterjee and Diaconis, can one simply find $p^*$ and obtain a sufficiently accurate approximation of the $\ergm$ by sampling $G(n,p^*)$? Unfortunately, no:
Theorem~\ref{thm:large_deviations} is in \emph{cut metric}, which turns out to be too weak to control total variation distance. Indeed, the paper \cite{bresler2018optimal} shows that $\tv(\ergm, G(n,p^*))\to 1$ even for the simple $\ergm$ consisting of edges and 2-stars.  
The approximation results of Eldan and Gross \cite{eldan2017exponential} show that even in the low-temperature regime, the $\ergm$ can be approximated in a certain Wasserstein metric by an appropriate mixture of stochastic block models, endowing it with richer structure compared to $G(n,p^{*})$. However, these results also do not imply approximation in the total variation distance.

The $\ergm$ is closely related to the ferromagnetic Ising model and, in fact, the $\ergm$ with the 2-star can be written as an Ising model. Just as for the $\ergm$, the Glauber dynamics is known to mix exponentially slowly for the Ising model at low temperatures. Nevertheless, there are Markov chains mixing in polynomial-time which can sample efficiently from arbitrary ferromagnetic Ising models based on random cluster dynamics \cite{guo2017random} and the Swendsen-Wang dynamics \cite{swendsen1987nonuniversal,ullrich2014swendsen}.
For the $\ergm$ it is not at all clear how to write down the corresponding random cluster model such that the associated random cluster dynamics mixes rapidly.
We instead pursue a more direct approach.

The starting point of our approach is the observation that if our aim is only to produce a sample from nearly the correct distribution, then there is no need for the dynamics to transition between all metastable states. 
 In order to implement this intuition, it is necessary to slightly modify the standard definition of mixing time of a Markov chain. The usual definition measures the distance to stationarity starting from a worst-possible initial state. Instead, we use the following definition.
\begin{definition}\label{def:metastable_mixing}
Given a Markov transition kernel $P$ with stationary $\pi^*$ we start from some initial distribution $\pi_0$ and 
say 
that $P$ is $(\pi_0,\pi^{*},\tau,\delta)$-\textbf{mixing} if for every $t \geq \tau$
\begin{equation}
    \label{eq:approx_mixing}
    \tv(\pi_0 P^t,\pi^{*}) \leq \delta\,.
\end{equation}
\end{definition}
Note that an immediate consequence of the data processing inequality for total variation is that $P$ is $(\pi_0,\pi^{*},\tau,\delta)$-mixing if and only if $\tv(\pi_0 P^{\tau},\pi^{*}) \leq \delta$.

The role of $\delta_0$ merits some discussion.
Incorporating the starting distribution into the definition of mixing time invalidates one of the basic lemmas: it is no longer true that the total variation decreases exponentially fast once the mixing time is exceeded. The basic reason is that in a Markov chain with multiple metastable states requiring exponential time to leave, any initial error in probability assigned to the metastable states might persist for exponential time. Thus, one might think of $\delta_0$ as capturing this initial (possibly unavoidable) error. Gheissari and Sinclair's work \cite{gheissari2021low} on mixing in low-temperature Ising models also considers mixing up to a TV distance of $\delta_0$. 

\sloppy
Our main result, Theorem~\ref{thm:main_result} given in Section~\ref{sec:results}, shows that whenever $\delta \geq \delta_0 = \exp(-\Omega_{\beta}(n))$, with $\pi_0$ being a mixture of $G(n,p^{*})$ for some carefully chosen distribution $p^{*}$, the Glauber dynamics for $\ergm$ is $(\pi_0,\pi^{*},C_{\beta}n^2\log\frac{n}{\delta},\delta)$-mixing even in the low temperature regime. That is, as long as as the target TV distance is $\geq \delta_0$, then the mixing time of the Glauber dynamics is $O(n^2\log n)$. This gives a counterpoint to the criticism of these models in \cite{bhamidi2011mixing} based on the difficulty of sampling these models at low temperature. The following is a corollary of Theorem~\ref{thm:main_result}.

\begin{theorem}
\label{thm:main_result_simple}
Suppose that $\ps$ is the unique global maximizer of $L_\beta$ and moreover that $L_{\beta}$ has nonzero second derivative at $p^{*}$. Let $\pi_0 := G(n,p^{*})$. There exist positive constants $c_{\beta}, C_{\beta}$ and $n_0(\beta)$, such that if $n >  n_0(\beta)$,
then whenever $\delta \geq \exp(-c_{\beta}n)$, the Glauber dynamics for $\mu_\beta$ is $\b(\pi_0,\mu,C_{\beta}n^2\log({n^2}/{\delta}),\delta\b)$-mixing.
\end{theorem}

In fact, Theorem~\ref{thm:main_result} is richer and shows that even when there are multiple global maximizers in the low temperature regime, we can sample efficiently from the conditional distribution of being close to any of these maximizers. In Theorem~\ref{thm:local_metastability}, we establish a surprising richness which can be present in the $\ergm$ at low temperature. Even within a small cut-metric neighborhood of the constant $p^{*}$ graphon where the $\ergm$ measure concentrates (and looks very close to $G(n,p^{*})$) \cite{chatterjee2013estimating}, we establish the existence of metastable states for certain choices of $\beta$ whose total probability is $\exp(-\Theta(n))$. The Glauber dynamics takes an exponentially long time to escape this set of metastable states. 
In contrast, Bhamidi et al.~\cite{bhamidi2011mixing} constructed metastable states as sets of graphs similar to $G(n,p)$, where $p$ was a local maximizer of $L_{\beta}$, and these have a total probability of $\exp(-\Omega(n^2))$.

\subsection{Overview and Proof Sketch}
We now give a high level overview of the ideas behind the main results, which are stated in Section~\ref{sec:results}.

\paragraph{Sufficient Conditions for Path Coupling} 
The large deviations results in \cite{chatterjee2013estimating} stated here as Theorem~\ref{thm:large_deviations} show that a sample from the $\ergm$ is w.h.p. close in cut-metric to some constant graphon with value $p^{*}$. Sufficient conditions established in \cite{bhamidi2011mixing} for the path coupling argument to work requires much stronger control on small subgraph counts than provided by the cut-metric: the increase in homomorphism density of $G$ in $X$ formed by adding any edge $e$ to $X$ (denoted by $\delN e G X$) must be approximately ${2n^{-2}|E(G)|}(p^{*})^{|E(G)|-1}$ for every fixed subgraph $G$. We will show that $\delN e G X$ indeed concentrates close to this value with probability $1-\exp(-\Omega_{\beta}(n))$. The coupling argument showing how this statement implies our main theorem is contained in Sections~\ref{sec:pfOverview} and~\ref{sec:pfMain}.

\paragraph{Fixed Point Equations for Subgraph Concentration}
Theorem~\ref{thm:large_deviations} shows that $X\sim \mu$ is close to the constant graphon with value $p^{*}$.  
Section~\ref{sec:unif_conc} reduces the task of showing concentration of $\delN e G X$ to show additionally that:
(1) every node degree uniformly concentrates close to $p^{*}$ and (2)
that the number of common neighbors of any two vertices $u$ and $v$ is close $n(p^{*})^2$. 
We show these two properties as follows. First, concentration of degrees is established using the cavity method, discussed momentarily.
For the second property, we make use of the concentration of degrees to derive a fixed point equation for the common neighbor counts and use the concentration results for fixed point equations established in \cite[Theorem 1.5]{chatterjee2007stein}. We note that the concentration results given by \cite[Theorem 1.5]{chatterjee2007stein} in themselves do not seem to be sufficient to establish the concentration of $\delN e G X$ and the concentration of degrees as established by the cavity method is essential. 

\paragraph{Cavity Method for Degrees}
% In order to show the required uniform concentration of node degrees, w
We use the cavity method as developed in Section~\ref{sec:cavity} to first show that conditioned on the exponential random graph $X$ being close to the constant $p^{*}$ graphon, the normalized degree of \emph{every} vertex concentrates close to $p^{*}$. Graphon convergence can show that \emph{most} vertices have degree close to $p^{*}$ (see Lemma~\ref{lem:almost_degree_conc}). 
To obtain the uniform concentration, we look at the law of the the edges emanating from a single vertex (called the cavity) conditioned on the rest of the graph being close to the constant graphon $p^{*}$. We show that the ``mean field" generated by the rest of the graph forces the cavity vertex to have degree close to $p^{*}$ with high probability. This is established in Theorem~\ref{thm:cavity_conc} and Corollary~\ref{cor:degree_conc} and is our main technical innovation.

\subsection{Discussion and Future Work}
Fast mixing of a Markov chain can be used to establish concentration of measure, central limit theorems, and estimation of the partition function.  Concentration of measure, CLTs and approximation by $G(n,p)$ (\cite{ganguly2019sub,reinert2019approximating}) have been explored in the literature for high-temperature $\ergm$ models. It would be interesting to consider their extension to low-temperature $\ergm$ via the approximate mixing established in this work. Maximum likelihood estimation often involves estimation of the partition function. There are multiple works (\cite{jerrum1993polynomial,sinclair1989approximate,gelman1998simulating,geyer1992constrained,kou2006equi}) that efficiently approximate partition functions of a parametric family with a given parameter $\beta$ by efficiently generating samples from the distribution for every choice of the parameter $\beta^{\prime}$. We leave open the problem of estimating the partition function of the low-temperature $\ergm$.

\subsection{Related Work}

Polynomial or quasipolynomial time mixing from a well-chosen initial distribution has been explored for the mean field Ising model by Levin et al. \cite{levin2010glauber} and more recently for the Ising model on the lattice $\mathbb{Z}^d$ by Gheissari and Sinclair~\cite{gheissari2021low}. 
Lubetzky and Sly~\cite{lubetzky2021fast} consider mixing from specific initial conditions for the $1$-dimensional Ising model and identify initial states which allow faster mixing than the worst case by constant factors. The work \cite{lindsey2021ensemble} considers an idea similar to ours in isolating the `modes' of a probability distribution in order to aid sampling. Their approach consider multiple, coupled, random instantiations of the Markov chains, which are all allowed to interact as the evolve, whereas our work considers a single instance of such a Markov chain.

\input{background}

\input{mainresults}

\input{mixing.tex}

\input{mainproof}

\input{keyThmPf}

\input{technical.tex}

\input{cavity}

\input{metaStabPf}

\section*{Acknowledgements}
%\begin{acknowledgements}
% Part of this work was done when D.N. was a graduate student at MIT and E.N was an undergraduate student at MIT. 
G.B. and D.N. gratefully acknowledge the hospitality of the Simons Institute for Theoretical Computer Science during Fall 2020 and Fall 2021. 
This work was supported in part by NSF CAREER award CCF-1940205 and NSF Award DMS-2022448. 
%\end{acknowledgements}


\bibliographystyle{spmpsci}      % mathematics and physical sciences
\bibliography{references}   % name your BibTeX data base

% \begin{thebibliography}{}
% %
% % and use \bibitem to create references. Consult the Instructions
% % for authors for reference list style.
% % Format for Journal Reference
% Author, Article title, Journal, Volume, page numbers (year)
% Author, Book title, page numbers. Publisher, place (year)
% \end{thebibliography}
\appendix

\input{techLems}

\end{document}

%% file: tempmacros.tex
\usepackage{amsthm}

\newtheorem{theorem}{Theorem}[section]
\newtheorem{lemma}[theorem]{Lemma}

\newtheorem{corollary}[theorem]{Corollary}

\newtheorem{innercustomthm}{Corollary}
\newenvironment{customthm}[1]
  {\renewcommand\theinnercustomthm{#1}\innercustomthm}
  {\endinnercustomthm}

\theoremstyle{definition}
\newtheorem{definition}[theorem]{Definition}
\newtheorem{example}[theorem]{Example}
\newtheorem{remark}[theorem]{Remark}

%% file: background.tex
\section{Background and Notation}
\label{sec:definition}
This section contains the basic definitions and notation for graphons and the Glauber dynamics. 

\subsection{Graphon Theory}\label{s:graphon}
Our work relies heavily on the theory of graph limits and graphons; our notation follows \cite{chatterjee2013estimating}. Let $\mathcal{W}$ denote the space of symmetric measurable functions $f : [0,1]^2 \rightarrow [0,1]$, where the space $[0,1]^2$ is endowed with the uniform probability measure. For $f, g \in \mathcal{W}$, define their \emph{cut distance} to be
$$\delta_\square(f, g) = \sup_{S, T \subset [0, 1]}\left |\int_{S \times T} \big(f(x, y) - g(x,y)\big)dxdy\right |\,,$$
where the supremum is over Borel measurable sets $S,T$. Define the equivalence relation $\sim$ on $\mathcal{W}$ by $f \sim g$ iff there exists a measure preserving bijection $\sigma:[0,1] \rightarrow [0,1]$ such that $f(x, y) = g(\sigma x, \sigma y) := g_\sigma(x, y)$. Let $\Tilde{\mathcal{W}}$ be the quotient space with respect to this equivalence relation. For $f \in \mathcal{W}$, let $\Tilde{f}$ denote its orbit in $\mathcal{W}$. A metric $\delta_\square$ on $\Tilde{\mathcal{W}}$ can now be defined as
$$\delta_\square(\Tilde{f}, \Tilde{g}) = \inf_{\sigma} \delta_\square(f, g_\sigma)\,.$$
An important fact in the theory of graph limits is that $(\Tilde{\mathcal{W}}, \delta_\square)$ is a compact metric space.

For a graph $X$ with vertex set $[n]$, we can associate the function $f^X \in \mathcal{W}$ where $f^X(x, y) = \ind_{(\lfloor nx\rfloor, \lfloor ny \rfloor) \in E(X)} = X_{\lfloor nx\rfloor, \lfloor ny \rfloor}$. We define its corresponding graphon to be $\Tilde{X} = \Tilde{f}^X \in \Tilde{\mathcal{W}}.$ Note that under this mapping, vertex isomorphic graphs correspond to the same element of $\Tilde{\mathcal{W}}$. We will denote the graphon with the constant value $p \in [0,1]$ by $p\mathbf{1}$.

For a finite simple graph $H$ with vertex set $[k]$ for some $k \in \mathbb{N}$ and a graphon $h \in \tilde{W}$, we define the homomorphism density $t(H, h)$ as
  $$  t(H, h) = \int_{[0, 1]^k} \prod_{(i,j) \in E(H)}h(x_i,x_j)dx_1\dots  dx_k\,. $$
In particular, the subgraph counts $N_i(X)$ appearing in the Hamiltonian $\mathcal{H}_\beta(X)$ are defined as the homomorphism densities 
\begin{equation} 
N_i(X) := t(G_i, \tilde{f}^X)\,.
  \label{eq:hom_den} 
\end{equation}
When emphasizing a particular graph $G$ we will also use the notation $N_G(X)$.

The classical theory of graph limits is too coarse to understand convergence of Markov chains because the cut metric does not control degrees of individual vertices (or neighborhoods of two vertices), which can have a large impact on the evolution of the Glauber dynamics.
The following quantities will allow us to establish a fine-grained understanding of the measure $\mu$.  Given a graph $X$ with $n$ vertices, whenever $u \in [n]$, we define the normalized degree 
\begin{equation}\label{eq:norm_degree}
   p_u(X) := \frac{\text{degree of vertex } u }{n} = n\int_{0}^{1}\int^{\frac{u}{n}}_{\frac{u-1}{n}}f^{X}(x,y)dxdy  \,.
\end{equation}
Similarly, define the normalized wedge count for nodes $u,v$ by
\begin{align}
    p_{uv}(X) &:= \frac{\text{number of vertices to which both $u$ and $v$ have edges}}{n}\nonumber \\ &= n^2 \int_{0}^{1}\int_{\frac{u-1}{n}}^{\frac{u}{n}}\int_{\frac{v-1}{n}}^{\frac{v}{n}}f^{X}(x,z)f^{X}(y,z)dx dy dz \,.\label{eq:norm_degree_2}
\end{align}

\subsection{Glauber Dynamics for the ERGM}
\label{sec:Glauber}

As described in the introduction, at each step of the Glauber dynamics for $\mu$ a pair of vertices $e=\{u,v\}$ is chosen uniformly at random from the $n\choose 2$ possibilities and the variable $X_e$ indicating presence of edge $e$ is updated according to the conditional probability 
\begin{equation}\label{e:phie}
\phi_e(X_{\sim e}) := \mathbb{E}_{X\sim \mu}\left[X_e\bigr| X_{\sim e}\right]\,.
\end{equation}

It will be useful to express the update probability in terms of subgraph counts.
For any graph $G = (V,E)$ define
\begin{equation}
    \delN e G X = N_G(X^{+e}) - N_G(X^{-e})
\end{equation}
and let
$$r_G(X,e) := \left(\frac{n^2\delN e G X}{2|E|}\right)^{\tfrac{1}{|E|-1}}\,.$$ 
The update probability can be expressed as
$$\phi_e(X_{\sim e}) = \mathbb{ E}_{X\sim\mu}[X_e|X_{\sim e}]=\frac{\exp(2\beta_0 + \sum_{l=1}^{K}2\beta_l |E_l|r_{G_l}(X,e)^{|E_l|-1})}{1+\exp(2\beta_0 + \sum_{l=1}^{K}2\beta_l |E_l|r_{G_l}(X,e)^{|E_l|-1})}\,.$$
% Therefore, our techniques hinge on controlling the quantity $\delN e G X$ and hence $r_G(X,e)$ in order to understand the Glauber dynamics. 
%\gb{Is a factor 2 missing in front of $\beta_0$?}

It follows from the definitions that if $X\sim G(n,p)$, then $r_G(X,e)\approx p$ with high probability.
If, conversely, it were the case that $r_G(X,e)=p$ for all $G_1,\dots, G_K$, then the update probability takes the following simpler form
(with some abuse of notation) given by
$\phi_{\beta}:[0,1] \to [0,1]$ with 
\begin{equation}\label{e:phi}
    \phi_{\beta}(p) = \frac{\exp(\sum_{i=0}^{K}2\beta_i|E_i| p^{|E_i|-1})}{1+\exp(\sum_{i=0}^{K}2\beta_i|E_i| p^{|E_i|-1})}\,.
\end{equation}
Let $M_{\beta}$ denote the set of global maximizers of $L_\beta$. Let $U_{\beta} \subset M_{\beta}$ be the global maximizers where the second derivative of $L_{\beta}$ is nonzero.
It can be shown that $p^{*} \in U_{\beta}$ only if $p^{*} = \phi_{\beta}(p^{*})$ and $\phi^{\prime}_{\beta}(p^{*}) < 1$.
What this implies is that $\ps$ is a \emph{stable fixed point}: if the chain is started at $X$ where $r_G(X,e)\approx \ps$ for all $G \in {G_1,\dots, G_K}$, then that continues to hold for exponentially many steps. This was shown in Lemma 17 from \cite{bhamidi2011mixing} and is stated in our paper as Lemma~\ref{lem:non_exit}. 
%\gb{say something about relation to $L_\beta$?}

We will also need to consider the Glauber dynamics for distributions $\pi$ that assign zero probability to some graphs. 

\begin{definition}[Glauber Dynamics]
\label{def:glauber_gen}
Given any probability distribution $\pi$ over $\Omega = \{0,1\}^{{{n}\choose{2}}}$, we define the Glauber dynamics with respect to $\pi$ to be the discrete time Markov chain over $\Omega$ with single step transitions (to obtain $X^{\prime}$) as follows:
\begin{enumerate}
    \item Pick a coordinate $E\in\left[{{n}\choose{2}}\right]$ uniformly at random. 
    \item Given the coordinate $E$, define $X^{\oplus E}$ to be $X$ with edge $E$ flipped, and let
    \begin{equation}
        X^{\prime} = \begin{cases} X^{\oplus E} \text{ with probability } \frac{\pi(X^{\oplus E})}{\pi(X^{\oplus E})+ \pi(X)} &\text{ if } \pi(X) \neq 0 \\
        X^{\oplus E} \text{ with probability $1$} &\text{ if } \pi(X) = 0 \\
        X \text{ with probability } \frac{\pi(X)}{\pi(X^{\oplus E})+ \pi(X)}
        \end{cases}
    \end{equation}
    % Where $X^{\oplus E}$ denotes flipping coordinate $E$.
\end{enumerate}
% 
%\gb{Any reason to express in terms of new notation $X^{\oplus E}$ rather than just $+E$ and $-E$?}
 The Glauber dynamics with respect to $\pi$ is reversible for $\pi$ (\cite{levin2017markov}), so in particular $\pi$ is stationary. When $\pi(X) > 0$ for every $X$, this reduces to the definition given in Definition~\ref{def:glauber}
\end{definition}

\subsection{Notation}
For any simple graph $X$, we denote by $V(X)$ its vertex set and by $E(X)$ its edge set. Given an unordered pair $e = (u,v)$ for some $u,v \in V(X)$, we define $X_e = X_{uv} = X_{vu} = \ind(e\in E(X))$. Without loss of generality, we take $V(X) = [n]$, where $|V(X)| = n$, and identify the space of finite simple graphs on $n$ vertices with the space $\Omega := \{0,1\}^{{{n}\choose{2}}}$, where the coordinates are indexed by tuples $(u,v)$ for $u<v$, $u,v \in [n]$. 
Throughout, we will reserve $u,v,w$ to denote vertices of size $n$ random graphs for large $n$ and $i,j,l$ to denote vertices of fixed graphs like $G_0,G_1,\dots$ above.
% The graph $X$ is then identified with the element of $\Omega$ whose value at co-ordinate $(u,v)$ is $X_{uv}$ for every $u,v \in [n]$ such that $u < v$. Henceforth, we will use this space $\Omega$ interchangeably with the space of graphs.  

By $X_{\sim e}$ we denote the graph formed by all edges other than the edge $e$. 
% i.e, for $X \in \Omega$, $X_{\sim e}$ is the vector of all the co-ordinates other than the one indexed by $e$. 
Given $X \in \Omega$, define $X^{+e} \in \Omega$ (resp. $X^{-e} \in \Omega$) by $(X^{+e})_{\sim e} = X_{\sim e}$ (resp. $(X^{-e})_{\sim e} = X_{\sim e}$) and $(X^{+e})_e = 1$ (resp. $(X^{-e})_e = 0$) i.e., we add (resp. remove) edge $e$ to the graph $X$. 

We use the standard asymptotic notation $O(\,\cdot\,),\Omega(\,\cdot\,)$, and $\Theta(\,\cdot\,)$. For $x,y \in \mathbb{R}^{+}$ $y = O_{\gamma}(x)$, we mean $y \leq C_{\gamma} x$ for some constant $C_{\gamma}$ which depends only on $\gamma$ (and similarly for $\Omega_{\gamma}$ and $\Theta_{\gamma}$). In the statement of the results, expressions of the form $\epsilon < c(\gamma)$ mean ``$\epsilon$ smaller than a constant depending only on $\gamma$" and $n > n_0(\gamma)$ means ``$n$ larger than a constant depending only on $\gamma$".

We will occasionally use the function $I:[0,1]\to -[1/2,0]$ given by $I(p) := \frac{1}{2}p\log p+\frac{1}{2}(1-p)\log (1-p)$. This is just $-1/2$ times the binary entropy function.

%% file: mainresults.tex
\section{Main Results}
\label{sec:results}
Fix $\beta \in \mathbb{R}\times \left(\mathbb{R}^{+}\right)^K$, and recall that $M_{\beta}$ denotes the set of global maximizers of $L_\beta$.
Let $U_{\beta} \subset M_{\beta}$ be the global maximizers where the second derivative of $L_{\beta}$ is nonzero. Throughout, we will always take $p^{*}\in U_{\beta}$. When $|U_{\beta}| = |M_{\beta}| = 1$, our main result, stated below in Theorem~\ref{thm:main_result}, shows that the Glauber dynamics for $\ergm$ when initialized at the $G(n,p^{*})$ distribution rapidly approximately mixes as long as the target total variation distance $\delta_0 \geq \exp(- c_{\beta} n)$ (see Definition~\ref{def:metastable_mixing}). Note that $|U_{\beta}| = |M_{\beta}| = 1$ even in the low temperature regime for Lebesgue almost all $\beta$.

In the case that $|U_{\beta}|>1$, we show that Glauber dynamics with the same initialization as above can efficiently and approximately sample from the $\ergm$ conditioned on being close in cut metric to the constant $p^{*}$ graphon. Note that Theorem~\ref{thm:large_deviations} shows that with a very large probability, a sample from the $\ergm$ is close to the constant $p$ graphon for some $p \in U_{\beta}$. If the probability of being close in cut-metric to each $p \in U_{\beta}$ under the measure $\mu$ is known, then we can initialize the Glauber dynamics to the correct mixture of $(G(n,p))_{p \in U_{\beta}}$ and show that it mixes rapidly as long as the target total variation distance $\delta_0 \geq \exp(-c_{\beta}n)$. In this work, we do not consider the problem of estimating these mixture probabilities.

% \gb{needs a little more interpretation / discussion}

% \gb{This suggests that maybe we can somehow solve even the general case with multiple maximizers... can we get self-reducibility somehow in order to find marginals and then approximate partition fn?}
% \dn{If we have an estimate for the partition function, we could do it, but the TV error of $\epsilon$ would be obtained after $poly(1/\epsilon)$ samples even if self reducibility works out}
 For $\eta >0$ denote the $\eta$-ball in cut metric around $p^*$ by
 $$
 \Ball \ps \eta := \{X \in \Omega: \delta_{\square}(X,p^{*}\mathbf{1}) \leq \eta\}\,.
 $$ 
 Let $P_{\eta}$ be the kernel of the Glauber dynamics with respect to the measure $\mu(\,\cdot\,| \Ball \ps \eta)$ and let $P$ be the kernel of the Glauber dynamics with respect to the measure $\mu$.  %The first result stated below shows that whenever $\pi_0$, the initial distribution, is set as $G(n,p^{*})$, the Glauber dynamics $P_{\eta}$ converges to  $\mu(\spacedot|\Ball \ps{\eta})$ rapidly up to a TV distance of $\exp(-\Omega(n))$ and the trajectory of $P_{\eta}$ is nearly indistinguishable from the trajectory of $P$ for an exponentially long time.  

\begin{theorem}
\label{thm:main_result}
 Let $\pi_0 := G(n,p^{*})$ for any $\ps\in U_\beta$. Let $\bar{X}_0\sim \pi_0,\bar{X}_1,\bar{X}_2,\dots,$ evolve according to $P_\eta$ and $X_0 \sim \pi_0,X_1,X_2,\dots$ evolve according to $P$. There exist positive constants $\eta_0(\beta)$, $c_{\beta,\eta}, C_{\beta,\eta}$ and $n_0(\beta,\eta)$ such that whenever $\eta < \eta_0(\beta)$ and  $n >  n_0(\beta,\eta)$, the following hold:
\begin{enumerate}
    \item $(X_0,X_1, \dots,X_T)$ can be coupled with $(\bar{X}_0,\dots,\bar{X}_T)$ such that with probability at least $1-TC(\beta,\eta)\exp(-c_{\beta,\eta}n)$, we have
    $$(X_0,X_1, \dots,X_T) = (\bar{X}_0,\bar{X}_1, \dots,\bar{X}_T)\,.$$ 
    \item Whenever $\delta \geq \exp(-c_{\beta,\eta}n)$, $P_{\eta}$ is $(\pi_0,\mu(\cdot| \Ball \ps \eta),C_{\beta,\eta}n^2\log({n^2}/{\delta}),\delta)$-mixing.
    \item If $|U_{\beta}| = |M_{\beta}| = 1$, then whenever $\delta \geq \exp(-c_{\beta,\eta}n)$, $P$ is $(\pi_0,\mu,C_{\beta,\eta}n^2\log({n^2}/{\delta}),\delta)$-mixing.
\end{enumerate}
\end{theorem}

This theorem shows that we can achieve metastable mixing by disregarding a portion of the state space of probability $\exp(-c_{\beta,\eta}n)$ under the measure $\mu(\cdot| \Ball \ps \eta)$. One might wonder if this is necessary, and in particular whether it is possible to improve the second item due to $ \Ball \ps \eta$ being possibly well-connected. We next answer this question in the negative and gain insight into the structure of the $\ergm$ measure at low temperature.

The paper \cite{bhamidi2011mixing} constructs metastable states where the graph is close to $G(n,p)$ for some $p$ which is a local maximizer of $L_{\beta}$, from which any local Markov chain takes $\exp(\Omega(n))$ time to escape. The large deviations theory based results established in \cite{chatterjee2013estimating} show that when $p$ is not the global minimizer of $L_{\beta}$, then these metastable states collectively have mass $\exp(-\Omega(n^2))$. One might hypothesize that the metastable states can be fully characterized by the behavior of local maximizers of $L_{\beta}$ and the cut-metric neighborhoods, and moreover that they have total mass $\exp(-\Omega(n^2))$ . 

Perhaps surprisingly, it turns out instead that the low-temperature $\ergm$ landscape is remarkably intricate even within the neighborhood $ \Ball \ps \eta$ around the global optimizer $p^*$. We show by construction that this set can contain multiple metastable states which collectively have mass $\exp(-\Theta(n))$ and from which the Glauber dynamics takes $\exp(\Omega(n))$ time to escape. These states are close in cut-metric to the constant graphon $p^{*}$ and it follows that cut-metric based large deviations analysis cannot capture the intricacies of Markov chain mixing in the $\ergm$ at low-temperatures.

\begin{example}\label{ex:t_ergm}
Suppose $K = 1$ and let $G_1$ be the triangle graph (i.e, the $3$ clique). Let $\sigma(x) := e^x/(1+e^x)$. There exist parameters $\beta_0,\beta_1 \in \mathbb{R}\times \mathbb{R}^{+}$ and real numbers $p_1^*\neq p_2^*$ such that:
\begin{enumerate}
    \item $p_1^{*}$ and $p_2^{*}$ satisfy $p_i^{*} = \sigma(2\beta_0 + 6\beta_1 (p_i^{*})^2)$, $U_{\beta} = \{p_1^{*}\}$, and $p_2^{*}$ is a local maximizer of $L_{\beta}$;
    \item There exists $q^{*} \in [0,1]$, $q^*\notin\{ p_1^*,p_2^*\}$, such that $q ^{*}= \sigma(2\beta_0 + 6\beta_1 q^* p_1^{*})$;
    \item Taking $f(x) = \sigma(2\beta_0 + 6\beta_1 x^2) $ and $g(x) = \sigma(2\beta_0 + 6\beta_1 x p_1^{*})$, we have $f^{\prime}(p_1^{*}) < 1$ and $g^{\prime}(q^{*}) < 1$.
\end{enumerate}
We numerically check that the choice $\beta_0 = -1.8$ and $\beta_1 = 2$ has $p_1^*, p_2^*$, and $q^{*}$ satisfying the relations above. As shown next, this turns out to imply metastability.
\end{example}

\begin{theorem}\label{thm:local_metastability}
Consider Example~\ref{ex:t_ergm} given above.  Let $\eta > 0$ be any small enough constant. Let the initial state $X_0$ be such that $(X_0)_{1j} \sim \mathsf{Ber}(q^{*})$ and $(X_0)_{ij} \sim \mathsf{Ber}(p_1^{*})$ for $i,j \neq 1$ and $i< j$ are independently distributed. Suppose $X_0,X_1,\dots$ is the trajectory of the Glauber dynamics with respect to $\mu$ with $\beta$ as given in Example~\ref{ex:t_ergm}.  Define the set of graphs $\Omega_{q,p}(\eta)$ for $q,p \in [0,1]$ by
 $$\Omega_{q,p}(\eta) := \{X\in \Omega:  \delta_{\square}(X,p) \leq {\eta}/{2}\text{ and } |p_1(X)-q| \leq \eta \}\,.$$
Then the following hold:
\begin{enumerate}
    \item The set $\Omega_{q^*,p_1^*}(\eta)$ is metastable: There is a constant $\alpha>0$ such that $$\mathbb{P}\left(\cap_{t \leq \exp(\alpha n)}\{X_t \in \Omega_{q^{*},p_1^{*}}(\eta)\}\right) \geq 1-\exp(-\Omega_{\eta}(n)) \,.$$
    \item The set $\Omega_{q^{*},p_1^{*}}(\eta)$ has sizable probability: $$\mu(\Omega_{q^{*},p_1^{*}}(\eta)) = \exp(-\Theta_{\eta}(n))\,.$$
    \item Most of the mass lies in $\Omega_{p_1^{*},p_1^{*}}(\eta)$: $$\mu(\Omega_{p_1^{*},p_1^{*}}(\eta)) \geq 1- \exp(-\Omega_{\eta}(n))\,.$$
\end{enumerate}
%\gb{Do we need to say something about how big $\eta$ is?}
\end{theorem}

The theorem is proved in Section~\ref{sec:ThmMetastable}.

In order to see why the set of states described in the theorem above are metastable, consider the first step of Glauber dynamics taking $X_0$ to $X_1$. The number of triangles formed by including an edge $e = (1,j)$ is approximately $n q^{*}p_1^{*}$ (which is $6np_1^{*}q^{*}$ after counting re-labelings), i.e., $\delN e {G_1} {X_0}\approx \frac{6p^{*}_1q^{*}}{n^2}$. Thus, the Glauber dynamics updates this coordinate to $1$ with probability $\approx \sigma(2\beta_0+6\beta_1q^{*}p_1^{*}) = q^{*}$. Similarly, if an edge $e = (i,j)$ is to be updated with $i,j \neq 1$, then the number of triangles formed is $n(p_1^{*})^2$ (which is $6np_1^{*}q^{*}$ after counting re-labelings) i.e., $\delN e {G_1} {X_0} \approx \frac{6(p_1^{*})^2}{n^2}$, and the probability of setting this coordinate to $1$ is $\approx \sigma(2\beta_0+6\beta_1(p_1^{*})^2) = p_1^{*}$. Therefore, the Glauber dynamics update still makes $X_1$ look approximately like the initial distribution. Not only that, but this is a stable fixed point, which follows from the conditions $f^{\prime}(p_1^{*}) < 1$ and $g^{\prime}(q^{*}) < 1$.
%\gb{I changed $\phi$ to $\sigma$ in this para}
%\gb{Might be good to state in terms of $\delN e G$ -- I started this, pls update the '???'}

%% file: mixing.tex
\section{Showing \Stable~Mixing for Glauber Dynamics }

\label{sec:pfOverview}

\subsection{Couplings, Contraction, and Mixing}
Consider a Markov chain over the finite state space $\mathcal{X}$ and with transition kernel $P$. Let $d:\mathcal{X}\times\mathcal{X}\to \mathbb{R}^{+}$ be such that $\sup_{x,y \in \mathcal{X}} d(x,y) \leq \dmax$ and $\inf_{ x\neq y} d(x,y) \geq \dmin$. We will use the following lemma to establish \stable~mixing, proved in Appendix~\ref{sec:deferred2}.
% 
%Let the sequence $X_0,X_1,\dots$ and $Y_0,Y_1,\dots$ be two trajectories of the Markov chain, defined on a common probability space with an arbitrary initialization for $X_0$ and $Y_0$. Suppose there exists a subset $A \subseteq \mathcal{X}\times \mathcal{X}$ such that whenever $(x,y) \in A$, there exists a distribution $Q_{xy}$ over $\mathcal{X}\times \mathcal{X}$ which couples $P(x,\cdot)$ and $P(y,\cdot)$ such that whenever $(X^{\prime},Y^{\prime}) \sim Q_{xy}$, we have
%$$\mathbb{E}d(X^{\prime},Y^{\prime}) \leq (1-\gamma)d(x,y)\,.$$

\begin{lemma}\label{lem:approx_mix} 
Let $A \subseteq \mathcal{X}\times \mathcal{X}$ be such that for $(x,y) \in A$ there exists a $\gamma$-contractive coupling $Q_{xy}$ of $P(x,\cdot)$ and $P(y,\cdot)$, i.e. for $(X^{\prime},Y^{\prime}) \sim Q_{xy}$ we have
$$\mathbb{E}d(X^{\prime},Y^{\prime}) \leq (1-\gamma)d(x,y)\,.$$
Then, given any jointly distributed $(X_0,Y_0) \in \mathcal{X}\times\mathcal{X}$, there exists a coupling between the trajectories $(X_k)_{k\geq 0}$ and $(Y_k)_{k\geq 0}$ of the Markov chain $P$ such that
$$\mathbb{E}d(X_{k+1},Y_{k+1}) \leq (1-\gamma)\mathbb{E}d(X_k,Y_k) + \dmax p_k\,,$$
where $p_k := \mathbb{P}((X_k,Y_k) \in A^{\complement})$. Unrolling this recursion, we conclude that
$$\mathbb{E}d(X_{K},Y_K) \leq \dmax\Big[(1-\gamma)^{K} + \frac{\sup_{k\leq K} p_k}{\gamma}\Big]\,.$$
%We conclude that $P$ is $(\pi_0,\pi^{*},K,\delta_K)$ approximately mixing with $\delta_K = \frac{\dmax}{\dmin}\left[(1-\eta)^{K} +\frac{\sup_{k\leq K}q_k+p_k}{\eta}\right]$
\end{lemma}

The following corollary is immediate from the coupling characterization of total variation.
\begin{corollary}
In the setting of Lemma~\ref{lem:approx_mix}, if additionally $\inf_{x\neq y}d(x,y)\geq \dmin$, then $\tv(X_K,Y_K)\leq \frac{\dmax}{\dmin}\b[(1-\gamma)^{K} + \gamma\inv \sup_{k\leq K} p_k\b]$.
\end{corollary}
In essence, the result above shows that whenever two trajectories can be coupled such that with high probability they lie in a set $A$ where a contractive coupling exists, then the laws of their iterates converge until a certain lower threshold. In particular, taking $Y_0$ to be drawn from the stationary distribution of $P$, we can establish \stable~mixing for $X_0,X_1,\dots$.

We will use the monotone coupling, defined next. 

\begin{definition}[Monotone coupling]\label{def:mon}When $P$ is the kernel of the Glauber dynamics with respect to $\mu$, the following coupling between $P(x,\cdot)$ and $P(y,\cdot)$ is called the \textbf{monotone coupling}. For any two $x,y \in \Omega$, we obtain the one step Glauber dynamics updates $X^{\prime},Y^{\prime}$ as follows:
\begin{enumerate}
    \item Pick the update edge $E \in {{[n]}\choose{2}}$ uniformly at random to be the same for both $X^{\prime}$ and $Y^{\prime}$.
    \item Draw $U \sim\unif( [0,1])$ independent of everything else and set
    \begin{equation}
        X^{\prime} = \begin{cases} x^{+E} \text{ if } U \in [0,\phi_{E}(x_{\sim E}))\\
         x^{-E} \text{ otherwise}
        \end{cases}
    % \end{equation}
 \qquad    \text{and}\qquad
        % \begin{equation}
        y^{\prime} = \begin{cases} y^{+E} \text{ if } U \in [0,\phi_{E}(y_{\sim E}))\\
         y^{-E} \text{ otherwise.}
        \end{cases}
    \end{equation}
\end{enumerate}
\end{definition}

For any two graphs $X,Y \in\Omega$, the relation $X \preceq Y$ denotes that $X_e \leq Y_e$ for every  $e \in {{[n]}\choose{2}}$.
It follows immediately from the definition of the monotone coupling that given $X\preceq Y$, if $X', Y'$ are obtained via the monotone coupling, then $X'\preceq Y'$ almost surely. We next identify a region of the state space over which the coupling is contractive (as required by Lemma~\ref{lem:approx_mix}).

\subsection{Control of Subgraph Counts Implies Contraction}

We will now follow the results established in \cite{bhamidi2011mixing} to show the path coupling of Glauber dynamics and use the notations they introduced. 
% To this end, for any graph $G = (V,E)$, we define $$r_G(X,e) := \left(\frac{\delN e G X}{2|E|n^{|V|-2}}\right)^{\tfrac{1}{|E|-1}}\,.$$ 
% 
Recall that the update probability under the Glauber dynamics for $\mu$ is given by
$$\phi_e(X_{\sim e}) = \frac{\exp(\beta_0 + \sum_{l=1}^{K}2\beta_l |E_l|r_{G_l}(X,e)^{|E_l|-1})}{1+\exp(\beta_0 + \sum_{l=1}^{K}2\beta_l |E_l|r_{G_l}(X,e)^{|E_l|-1})}\,,$$ 
where
$$
r_G(X,e) := \Big(\frac{n^2 \delN e G X}{2|E(G)|}\Big)^{\tfrac{1}{|E(G)|-1}}$$
and $\delN e G X = N_G(X^{+e}) - N_G(X^{-e})$.

Let $\mathbb{G}_L$ denote the set of finite simple graphs with at most $L$ vertices (omitting the graph with $1$ edge and $2$ vertices), 
% where $L$ is a fixed constant that does not depend on $n$. For the rest of this work, we will take 
where $L$ is a fixed constant satisfying $L > \max_{i\leq K}|V_i|$. 
% Define 
% \begin{equation*}
%     r_{\max}(X) := \sup_{\substack{e \in {{n}\choose{2}} \\ G \in \mathbb{G}_L}} r_G(X,e)
% % \end{equation}
% \qquad \text{and}\qquad
% % \begin{equation}
%     r_{\min}(X) := \inf_{\substack{e \in {{n}\choose{2}} \\ G \in \mathbb{G}_L}} r_G(X,e)\,.
% \end{equation*}
Define the set 
\begin{equation}\label{e:Omega0}
    \Gam \ps \eps := \b\{X: r_G(X,e) \in [p^{*}-\epsilon,p^{*}+\epsilon]\text{ for all } e\in \textstyle{[n]\choose 2} \text{ and } G\in  \mathbb{G}_L\b\}\,.
\end{equation}
Note that whenever $X \in \Gam \ps \eps$, $\phi_{e}(X_{\sim e}) \approx \phi_{\beta}(p^{*}) \approx p^{*}$. That is, each edge updates approximately like $G(n,p^{*})$. The significance of $\Gam \ps \eps$ is that in this set the monotone coupling is contractive, as shown in \cite[Lemma 18]{bhamidi2011mixing} and stated next. We will additionally state a theorem in the next subsection that $\Gam \ps \eps$ has high probability under $\mu(\cdot|\Ball \ps \eta)$.
% We also show that this set $\Gam \ps \eps$ has a high probability under the measure $\mu(\cdot|\Ball \ps \eta)$ in Theorem~\ref{thm:final_unif_conc}.

\begin{lemma}[Contraction within $\Gam \ps \eps$, \cite{bhamidi2011mixing}]%\dn{check this lemma alone}
\label{lem:mono_coup}
Let $p^{*}\in U_{\beta}$, $\epsilon > 0$ small enough as a function of $\beta$, and $n$ large enough as a function of $\beta,\epsilon$. Let $\mathcal{A} := \{ (x,y)\in \Gam \ps \eps\times\Gam \ps \eps : x\preceq y\}$. 
% Denote the monotone coupling between $P(X,\cdot)$ and $P(Y,\cdot)$ as $Q^{\mathsf{mon}}_{X,Y}$. 
Let $X',Y'$ be obtained from $x,y$ via one step of the Glauber dynamics under the monotone coupling.
There is a constant $c(\beta,\epsilon)>0$ such that if $(x,y)\in \mathcal{A}$, then
% $Q^{\mathsf{mon}}(X,Y)$ is contractive with respect to the Hamming distance with $\gamma = {c(\beta,\epsilon)}/{n^2}$ where $c(\beta,\epsilon) > 0$ whenever $(X,Y) \in \mathcal{A}$ (as defined above Lemma~\ref{lem:approx_mix}). 
$$
\mathbb{E}d_H(X^{\prime},Y^{\prime}) \leq \B(1-\frac{c(\beta,\epsilon)}{n^2}\B)d_H(x,y)\,,
$$
and moreover, 
% $(X^{\prime},Y^{\prime})\sim Q_{XY}$ satisfy 
$X^{\prime} \preceq Y^{\prime}$ almost surely. %whenever $(X,Y) \in \mathcal{A}$.
\end{lemma}

\subsection{Key Theorem}

We first recall Theorem~\ref{thm:large_deviations} which state that $\cup_{\ps \in U_{\beta}}\Ball \ps \eta$ has probability $1-\exp(-\Omega(n^2))$ under the measure $\mu$, that is, most of the mass of $\mu$ is concentrated in the cut-metric balls $\Ball \ps \eta$.
% , which allows us to concentrate within these balls to understand the measure $\mu$. 
The following theorem 
shows that $\mu\b(\cdot\b|\Ball \ps \eta\b)$ concentrates over the set $\Gam \ps \eps$, where path coupling is possible (as per Lemma~\ref{lem:mono_coup}).

\begin{theorem}\label{thm:final_unif_conc}
Suppose $p^{*} \in U_{\beta}$. Given $\epsilon > 0$, we can pick $\eta < c(\beta,\epsilon)$ such that
    $$\mu\b(\Gam \ps \eps\b|\Ball \ps \eta\b)\geq 1- C(\eta,\epsilon,\beta)\exp\big(-\Omega_{\beta,\epsilon,\eta}(n)\big)\,.$$
\end{theorem}

We prove the theorem in Section~\ref{sec:unif_conc}, modulo lemmas proved via the cavity method in Section~\ref{sec:cavity}. 

\begin{remark}
Notice that $r_G(X,e) \in [p^{*}-\epsilon,p^{*}+\epsilon]$ uniformly for every $e$ is not implied by $\delta_{\square}(\tilde{X},p^{*}) < \eta$ (for any constant $\eta>0$). An example is given in Theorem~\ref{thm:local_metastability}, where metastability occurs despite being close to the $p^{*}$ graphon with high probability: The edges emanating from a single vertex prevent uniform concentration of $r_G(X,e)$ in the set $[p^{*}-\epsilon,p^{*}+\epsilon]$, but the single vertex neighborhood has a vanishingly small impact on $\delta_{\square}(\tilde{X},p^{*})$. 
\end{remark}

While it can easily be proved directly, the following is also a corollary of the above theorem. 
\begin{corollary}
    Fix any $\eps>0$. Then there exists $\delta_0(\eps)>0$ such that for all $0<\delta<\delta_0(\eps)$, if $Z\sim G(n,\ps+\delta)$, then
    $
    \P (Z\in \Gam \ps \eps) \geq 1 - \exp(\Omega_{\eps}(n))
    $.
\end{corollary}
%\gb{this was added here, it wasn't ever explicitly stated, right?}

%\gb{say maybe not this -- which is exactly thm statement, but can we say how this will be used for us? i.e., we may do ?? instead of ??, or some such thing? I want to explain to reader how this is reducing our challenge. Now that this is shown, what does that do to make our life easier? 
%The related question is where are we getting $X_t$ from $\mu\b(\b|\Ball \ps \eta\b)$?}

%\gb{We hsould also emphasize that the $r_G$ are very detailed information that is much finer than the cut metric. So the surprise here is that we get a strong control.}
% \begin{theorem}\label{thm:final_unif_conc}
% Suppose $p^{*} \in U_{\beta}$. Given $\epsilon > 0$, we can pick $\eta < c(\beta,\epsilon)$ and for every graph $G = ([k],E)$ a constant $C_{G,\beta}>0$ such that
%     \begin{align}\mu\Big(\sup_{G} \frac{1}{C_{G,\beta}}\big|n^2\delN e G X - 2|E(G)|(p^{*})^{|E(G)|-1}\big| >\epsilon \,\Big| \Ball \ps\eta\Big) \nonumber
%     % \\&\qquad\qquad
%     \leq C(\eta,\epsilon,\beta)\exp\big(-\Omega_{\beta,\epsilon,\eta}(n)\big)\,.\end{align}
% \end{theorem}
% \gb{maybe emphasize that the content of the theorem is the fact that this holds for all $n$}

\subsection{Metastability}

 We intend to invoke Lemma~\ref{lem:approx_mix} to show approximate mixing and prove Theorem~\ref{thm:main_result}. The prior subsection shows that  $G(n,p^{*})$ and at $\mu\b(\cdot\b|\Ball \ps \eta\b)$ are both within the set $\Gam \ps \eps$ with high probability. We show that the Glauber chains with these initializations do not leave $\Gam \ps \eps$ with probability $1-\exp(-\Omega(n))$ until time $\exp(\Omega(n))$. Some intuition behind this was given in Section~\ref{sec:Glauber}.

We next state Lemma 17 from \cite{bhamidi2011mixing}, after adapting it to our situation. 

\begin{lemma}[%Modification of Lemma 17 in \cite{bhamidi2011mixing}
Staying in $\Gam \ps \eps$]\label{lem:non_exit}
Let $\epsilon > 0$ be a small enough constant independent of $n$ and suppose $p^{*} \in U_{\beta}$. Let $X_0,X_1,\dots$ evolve according Glauber dynamics with respect to the measure $\mu$. 
% For some large enough $L \in \mathbb{N}$ \gb{Can we just take $L$ as above?}, independent of $n$, 
If 
% $X_0$ is such that $p^{*}-\epsilon \leq r_{\min}(X_0) \leq r_{\max}(X_0) \leq p^{*} + \epsilon$, 
$X_0\in \Gam \ps \eps$,
then for some $\alpha=\alpha(\beta,\epsilon,L)$, we have
$$
\P\b(X_t\in \Gam \ps {2\eps}\text{ for all }t\leq e^{\alpha n}\b) \geq 1- \exp\big(-\Omega_{\beta,L,\epsilon}(n)\big)\,.
$$
% $$\mathbb{P}\Big(\sup_{t \leq e^{\alpha n}} r_{\max}(X_t) \geq  p^{*} + 2\epsilon\Big) \leq \exp\big(-\Omega_{\beta,L,\epsilon}(n)\big)$$
% and
% $$\mathbb{P}\Big(\inf_{t \leq e^{\alpha n}} r_{\min}(X_t) \leq  p^{*} - 2\epsilon\Big) \leq \exp\big(-\Omega_{\beta,L,\epsilon}(n)\big) \,.$$
\end{lemma}

% Letting $P$ be the kernel of the Glauber dynamics with respect to $\mu$, we will call this coupling between $P(X,\cdot)$ and $P(Y,\cdot)$ the \textbf{monotone coupling} and denote the joint law of $(X^{\prime},Y^{\prime})$ by \gb{?}

% \begin{equation}\label{e:unifSubG}
%  \usc=\B\{ \sup_{G} \frac{1}{C_{G,\beta}}\big|n^2\delN e G X - 2|E(G)|(p^{*})^{|E(G)|-1}\big| \leq \epsilon  \B\}
% \end{equation}

% \begin{lemma}
%     $\Gam \ps \eps\subseteq  \usc$ \gb{Immediate, by definition}
% \end{lemma}

%\gb{This is a lot: ``for some $\delta >0$, small enough as a function of $\beta,\epsilon$, and $\eta$ large enough as a function of $\delta,\epsilon$ and $\beta$, and $n$ large enough as a function of $\delta,\beta,\eta,\epsilon$," perhaps we can just introduce some simple notation for it and have it just appear once. }

The proofs of the following lemmas are given in Appendix~\ref{sec:deferred2}.

\begin{lemma}[$G(n,\ps\pm \eps)$ sandwich]\label{lem:stoc_dom}
Let $p^{*}\in U_{\beta}$. Let constants $\epsilon,\eta > 0$ be such that $\epsilon < \eps_0(\beta)$, $\eta < \eta_0(\beta,\epsilon)$, and $n > n_0(\beta,\epsilon,\eta)$. Let $X \sim \mu(\cdot|\Ball \ps \eta)$, $\bar{Y} \sim G(n,p^{*}+\epsilon)$, and $\underline{Y} \sim G(n,p^{*}-\epsilon)$. Then, there exists a coupling between $X,\bar{Y}$, and $\underline{Y}$ such that with probability at least $1-\exp(-\Omega_{\beta,\eta,\epsilon}(n))$
$$
\underline{Y}\preceq X\preceq \bar{Y}\,.
$$
%\gb{I put this display above into one line, instead of two separate comparisons, fine right?}
% Similarly, there exists a coupling between $\underline{Y} \sim G(n,p^{*}-\epsilon)$ and $X$ such that with probability at least $1-\exp(-\Omega_{\beta,\eta,\epsilon}(n))$
% $$ X\succeq \underline{Y}\,.$$
\end{lemma}

\begin{lemma}[Staying in $\Ball \ps {\eta/{2}}$]\label{lem:trajectory_inside_ball}
    Suppose $p^{*} \in U_{\beta}$, $\eta > 0$ such that $\eta < \eta_0(\beta)$ and $n > n_0(\eta,\beta)$. Let $X_0 \sim G(n,p^{*})$ and generate the trajectory $X_0,\dots,X_T$ via Glauber dynamics with respect to $\mu$. The entire trajectory $X_0,\dots,X_T$ stays within the ball $\Ball \ps {\eta/{2}}$ with probability at least $1-TC(\beta,\eta)\exp(-c(\beta,\eta)n)$.
\end{lemma}

%% file: mainproof.tex
\section{Proof of Main Result, Theorem~\ref{thm:main_result}}
\label{sec:pfMain}
% Before we begin the proof, we state the following theorem which is essential for constructing a path coupling. 

% It will be proved in Section~\ref{sec:unif_conc} after establishing the cavity method in Section~\ref{sec:cavity}. 

% \begin{theorem}\label{thm:final_unif_conc}
% Suppose $p^{*} \in U_{\beta}$. Given $\epsilon > 0$, we can pick $\eta < c(\beta,\epsilon)$ \gb{is this for all $\eta$? where are the quantifiers?} and for every graph $G = ([k],E)$ a constant $C_{G,\beta}>0$ such that
%     \begin{align}\mu\Big(\sup_{G} \frac{1}{C_{G,\beta}}\big|n^2\delN e G X - 2|E(G)|(p^{*})^{|E(G)|-1}\big| >\epsilon \,\Big| B(\emptyset,p^{*},\eta)\Big) \nonumber
%     % \\&\qquad\qquad
%     \leq C(\eta,\epsilon,\beta)\exp\big(-\Omega_{\beta,\epsilon,\eta}(n)\big)\,.\end{align}
% \end{theorem}
% \gb{maybe emphasize that the content of the theorem is the fact that this holds for all $n$}

We now show how the main theorem follows from the various lemmas stated in the last section.
Recall that $\pi_0 := G(n,p^{*})$, $\bar{X}_0\sim \pi_0,\bar{X}_1,\bar{X}_2,\dots,$ is a trajectory of the Markov chain $P_\eta$ and $X_0 \sim \pi_0,X_1,X_2,\dots$ is a trajectory of the Markov chain $P$.

\subsection{Proof of Theorem~\ref{thm:main_result}, Part 1}

    % \item Borrowing notation from Lemma~\ref{lem:non_exit} with some arbitrary but fixed $L$, we first consider $X_0 \sim \pi_0$. We know for $G(n,p^{*})$ and any $\epsilon >0$ small enough (\dn{cite references here}) that $\mathbb{P}(r_{\max}(X_{0}) \geq p^{*} + \epsilon ) \leq \exp(-\Omega_{\beta,L,\epsilon}(n))$ 
    % and 
    % $\mathbb{P}(r_{\min}(X_{0}) \leq p^{*} - \epsilon ) \leq \exp(-\Omega_{\beta,L,\epsilon}(n))$. 
    
    % Therefore, by an application of Lemma~\ref{lem:non_exit}, we conclude that with probability at-least $1-T\exp(-\Omega_{\beta,\epsilon,L}(n))$, we must have:
    % $$p^{*}-2\epsilon \leq \inf_{ 0 \leq t \leq T}r_{\min}(X_{t}) \leq \sup_{ 0 \leq t \leq T}r_{\max}(X_{t})\leq p^{*}+2\epsilon$$
    % \dn{need to check that $r_{\max}$ and $r_{\min}$ are indeed sufficient for the inverse function lemma}
    % Now, using the inverse counting lemma (Lemma~\ref{lem:inv_count}), we conclude by taking $L$ large enough and $\epsilon$ small enough that, with probability at-least $1-T\exp(-\Omega_{\beta,\eta}(n))$, we must have:
    
    % $$X_0,\dots,X_T \in B(\emptyset,p^{*},\frac{\eta}{2})\,.$$
    % \dn{check these carefully}
    % Under this event, whenever $n$ is large enough compared to $\eta$, this means that the updates of Glauber dynamics with respect to $\mu$ and $\mu(|B(\phi,p^{*},\eta)
    % )$ are exactly the same, which gives us the coupling \dn{elaborate more on this}
    We will couple the trajectories $\bar{X}_0,\dots,\bar{X}_T$ and $X_0,\dots,X_T$ such that the event $E:= \{(\bar{X}_0,\bar{X}_1,\dots,\bar{X}_T) \neq (X_0,\dots,X_T)\}$ satisfies $E \subseteq \cup_{t=0}^{T}\b\{X_t \in \b(\Ball \ps {\eta/{2}}\b)^{\complement}\b\}$. We can then conclude the result from Lemma~\ref{lem:trajectory_inside_ball}. The main observation is that whenever $n$ is large enough as a function of $\eta$, if $X \in \Ball \ps {\eta/{2}}$, then $P_{\eta}(X,\cdot) = P(X,\cdot)$. We construct the following coupling:
    \begin{enumerate}
        \item $X_0 = \bar{X}_0$ almost surely.
        \item $X_{t+1}, \bar{X}_{t+1}$ are drawn from the TV optimal coupling between $P(X_t,\cdot)$ and $P_{\eta}(\bar{X}_t,\cdot)$.
    \end{enumerate}
    It is clear that $\{X_{t+1}\neq \bar{X}_{t+1}\} \subseteq \{X_t \neq \bar{X}_t\}\cup\b\{X_t \in \b(\Ball \ps {\eta/{2}}\b)^{\complement}\b\}$. Now, noting that $\{X_0 \neq \bar{X}_0\}$ is the empty event, we conclude that $\{X_{1} \neq \bar{X}_1\} \subseteq \bigr\{X_0 \in \b(\Ball \ps {\eta/{2}}\b)^{\complement}\bigr\}$. An induction argument with the same basic step shows that $E \subseteq \cup_{t=0}^{T}\bigr\{X_t \in \b(\Ball \ps {\eta/{2}}\b)^{\complement}\bigr\}$.
\qed

\subsection{Proof of Theorem~\ref{thm:main_result}, Part 2}
     Let $\bar{Y}_0 \sim \mu\left(\cdot\bigr|\Ball \ps \eta\right)$ and consider the trajectory $\bar{Y}_0,\bar{Y}_1,\dots,\bar{Y}_T$ with respect to the transition kernel $P_{\eta}$.
      Similarly, let $Y_0 \sim \mu\left(\cdot\bigr|\Ball \ps \eta\right)$, but with the trajectory $Y_1,\dots, Y_T$ generated with respect to the transition kernel $P$. Using Lemma~\ref{lem:outer_prob} to bound $\mathbb{P}(\bar{Y}_t \in \Ball \ps \eta)$, a similar proof as in Item 1 shows that
    \begin{equation}\label{eq:unclean_coup_1}
    \tv\b((\bar{Y}_0,\dots,\bar{Y}_T), (Y_0,\dots,Y_T)\b) \leq T \exp\b(-\Omega_{\beta,\eta}(n^2)\b)\,.
    \end{equation}
    From Item 1, we have
    \begin{equation}\label{eq:unclean_coup_2}\tv((\bar{X}_0,\dots,\bar{X}_T), (X_0,\dots,X_T)) \leq T \exp(-\Omega_{\beta,\eta}(n))\,.
    \end{equation}
    
    These last two displays allow us to consider the total variation distance between the distributions of $X_T$ and $Y_T$ instead of $\bar{X}_T$ and $\bar{Y}_T$. Let $\epsilon > 0$ be small enough to satisfy the conditions in Lemma~\ref{lem:non_exit}. By Lemma~\ref{lem:stoc_dom} and Theorem~\ref{thm:final_unif_conc}, we conclude that for some $0<\delta <\delta_0(\beta,\epsilon,\eta)$, $\eta<\eta_0(\delta,\epsilon,\beta)$, and $n>n_0(\delta,\beta,\eta,\epsilon)$, 
    we can couple $Z_0 \sim G(n,p^{*}+\delta)$ with $X_0$ and $Y_0$ such that with probability at least $1-\exp(-\Omega_{\beta,\delta,\eta,\epsilon}(n))$, the following hold:
    \begin{enumerate}
        \item     $X_0 \preceq Z_0$ and $Y_0 \preceq Z_0$, and
        % \item $\{r_{\max}(X_0),r_{\min}(X_0),r_{\max}(Y_0),r_{\min}(Y_0),r_{\max}(Z_0),r_{\min}(Z_0)\} \subset [p^{*}-\epsilon,p^{*}+\epsilon]$.
                \item $X_0,Y_0, Z_0\in\Gam \ps \eps $\,.
    \end{enumerate}
    % The second item for $X_0, Z_0$ follows from the proof of  Lemma~\ref{lem:trajectory_inside_ball} and for $Y_0$ by Lemma~\ref{lem:non_exit}.

    % \gb{Actually for all three by the main theorem no? Neither of these lemmas have much to do with this, right?}

Now, we consider the Markov chain $Z_0,Z_1,\dots,Z_T$ with respect to $P$. We consider the monotone coupling between $Y_0,\dots,Y_T$ and $Z_0,\dots,Z_T$ as in Definition~\ref{def:mon}: conditional on $Y_0 \preceq Z_0$, we have $Y_t \preceq Z_t$ almost surely for every $t \leq T$. Recall $\Gam \ps \eps$ from \eqref{e:Omega0} and
the set $\mathcal{A} := \{ (x,y)\in \Gam \ps \eps\times\Gam \ps \eps : x\preceq y\}$ defined in Lemma~\ref{lem:mono_coup}. Now, $(Y_t,Z_t) \in \mathcal{A}^{\complement}$ implies either $Y_0 \not\preceq Z_0$, $Y_t \not\in \Gam \ps \eps$, or $Z_t \not\in \Gam \ps \eps$.  The first of these is ruled out by monotonicity of the coupling. For the latter two, Lemma~\ref{lem:non_exit} shows that starting in $\Gam \ps \eps$ the trajectory stays there for some time, and it follows that
$$\mathbb{P}\big((Y_t,Z_t)\in \mathcal{A}^{\complement}\big) \leq \exp\big(-\Omega_{\beta,\delta,\epsilon,\eta}(n)\big)\,.$$

%\gb{Confusing having $\gamma$ below -- contraction lemma is stated in terms of constant}

Applying Lemmas~\ref{lem:approx_mix} and~\ref{lem:mono_coup} with Hamming distance $d_H$ over $\Omega$, we conclude that whenever $T \leq \exp(c_0 n)$ for small enough $c_0$ as a function of $\epsilon,\delta,\eta,\beta$
$$\mathbb{E}d(Y_T,Z_T) \leq \frac{n^2}{2}\left[\Big(1-\frac{c(\epsilon,\beta)}{n^2}\Big)^{T} + \frac{n^2}{c(\epsilon,\beta)}\exp(-\Omega_{\eta,\beta,\epsilon,\delta}(n))\right]\,.$$
Similarly, we have
$$\mathbb{E}d(X_T,Z_T) \leq \frac{n^2}{2}\left[\Big(1-\frac{c(\epsilon,\beta)}{n^2}\Big)^{T} + \frac{n^2}{c(\epsilon,\beta)}\exp(-\Omega_{\eta,\beta,\epsilon,\delta}(n))\right]\,.$$
Combining the two displays above, the coupling characterization of TV distance implies 
\begin{align}
    \tv(X_T,Y_T) &\leq \mathbb{P}(X_T\neq Y_T) \leq \mathbb{P}(X_T\neq Z_T) + \mathbb{P}(Z_T\neq Y_T)\nonumber \\
    &\leq \mathbb{P}(d(X_T,Z_T) > 1) + \mathbb{P}(d(Y_T,Z_T) > 1) \nonumber \\
    &\leq \mathbb{E}d(X_T,Z_T) + \mathbb{E}d(Y_T,Z_T) \nonumber\\
    &\leq n^2\left[\Big(1-\frac{c(\epsilon,\beta)}{n^2}\Big)^{T} + \frac{n^2}{c(\epsilon,\beta)}\exp\big(-\Omega_{\eta,\beta,\epsilon,\delta}(n)\big)\right]\label{eq:clean_coup}\,.
\end{align}

Now, we will allow $\epsilon,\delta$ to be small enough constants as a function of $\beta$ such that Lemmas~\ref{lem:non_exit} and~\ref{lem:mono_coup} hold.  Whenever $\eta$ is small enough as a function $\beta$ and $n$ is large enough, combining Equations~\eqref{eq:unclean_coup_1},~\eqref{eq:unclean_coup_2}, and~\eqref{eq:clean_coup}, yields
$$\tv(\bar{X}_T,\bar{Y}_T) \leq n^2\left[\Big(1-\frac{c(\epsilon,\beta)}{n^2}\Big)^{T} + \frac{n^2T}{c(\epsilon,\beta)}\exp\big(-\Omega_{\eta,\beta}(n)\big)\right]\,.$$
This yields the second part of the theorem statement by considering $T = C_{\beta,\eta}n^2 \log(n^2/\delta)$.
%along with Lemma~\ref{lem:nec_suff_am}. 
  \qed 
  
%  \gb{made various edits in this proof of 2nd part, check that still 'compiles'}
   
   \subsection{Proof of Theorem~\ref{thm:main_result}, Part 3}
  Whenever $|U_{\beta}| = 1$, the concentration result of \cite{chatterjee2013estimating} (stated here as Theorem~\ref{thm:large_deviations}) implies that $$\tv\big(\mu,\mu\b(\cdot|\Ball \ps \eta\b)\big) \leq \exp\big(-\Omega_{\beta,\eta}(n^2)\big)\,.$$
   Meta-stable mixing to $\mu$ follows from the second part and the triangle inequality.  \qed

%% file: keyThmPf.tex
\section{Uniform Subgraph Concentration and Proof of Theorem~\ref{thm:final_unif_conc}}
\label{sec:unif_conc}

% Approximation of $\ergm$ with respect to cut metric is too coarse-grained to prove mixing of Markov chains. This is substantiated by Theorem~\ref{thm:local_metastability}, which constructs meta-stable states that are arbitrarily close to the constant graphons $p^{*} \in U_{\beta}$, yet have exponentially small probability under the $\ergm$ distribution.

% To show fast mixing for the restricted distribution $\mu(\spacedot|\Ball \ps \eta)$, the key ingredient 
% is Theorem~\ref{thm:final_unif_conc}, which states that 
%  with high probability over $X \sim \mu(\spacedot|\Ball \ps \eta)$ the quantity $r_G(X,e)$ is close to $\ps$ for every $e$ and small graph $G$.
% Equivalently,

In this section we reduce the proof of Theorem~\ref{thm:final_unif_conc} on the concentration of 
$\delN e G X:= N_G(X^{+e}) - N_G(X^{-e})$ to control of both the vertex degrees $p_u(X)$ and common neighbors $p_{uv}(X)$. The latter are stated here as corollaries and will be proved via the cavity method in Section~\ref{sec:cavity}.

Theorem~\ref{thm:final_unif_conc} states that if we sample $X \sim \mu\b(\cdot\b|\Ball \ps \eta\b)$, then
$X\in \Gam \ps \eps$ with probability $1- C(\beta,\epsilon,\eta)\exp\big(-\Omega_{\beta,\epsilon,\eta}(n)\big)$.
Unpacking the definitions, it suffices to show that for some small enough $h(\epsilon,L)$, 
\begin{equation}\label{eq:delt_conc}
    \B|\delN e G X - \frac{2|E(G)| (p^{*})^{|E(G)|-1}|}{n^2}\B| \leq \frac{h(\epsilon,L)}{n^2}
\end{equation} for every $G \in \mathbb{G}_L$ and $e \in {{[n]}\choose{2}}$, with probability at least $1- C(\beta,\eta,\epsilon,L)\exp(-\Omega_{\beta,\eta,\epsilon,L}(n))$.
% , since then $|r_G(X,e) - p^{*}| <\epsilon$ for every edge $e$ and $G \in \mathbb{G}_L$.

% $\delN e G X:= N_G(X^{+e}) - N_G(X^{-e})$, the change in the homomorphism density when the edge $e$ is flipped, is approximated as \begin{equation}\label{e:Napprox}
%  \delN e G X \approx   {2|E(G)|(p^{*})^{|E(G)|-1}}/{n^2}
% \end{equation} 
% for some $p^{*} \in U_{\beta}$, every edge $e$, and all small enough $G$. 
% % However, it is not clear that such a set of states $X$ has high probability under the measure $\mu$, especially in the low temperature regime. 
% The concentration established for $\delta_{\square}(\tilde{X},p^{*})$ in Theorem~\ref{thm:large_deviations} is too weak to show that \eqref{e:Napprox} holds uniformly for every edge $e$. 

% \gb{make eqn above quantitative so that can just proceed with precise bounds}

% The goal of this section is to translate control of $\delN e G X$ to control of degrees and then to outline the steps for achieving the latter. 
We start with a lemma (proved in Appendix~\ref{sec:deferred_proofs}) which shows that \eqref{eq:delt_conc} (and hence Theorem~\ref{thm:final_unif_conc}) follows from uniform control of both the vertex degrees $p_u(X)$ and common neighbors $p_{uv}(X)$.
Some notation is needed.
Given a fixed graph $G = ([k],E)$ and vertices $i,j \in [k]$ such that $(i,j)\in E$, let
\begin{equation}\label{e:pairDeg}
    E_{ij}(G) := \b\{l : l \in [k], l\neq i, l\neq j, (l,i) \in E, (l,j)\in E\b\}\quad \text{and} \quad d_{ij}(G) = |E_{ij}(G)|\,.
\end{equation} In words, $d_{ij}(G)$ is the number of common neighbors of vertices $i$ and $j$ in $G$.

\begin{lemma}\label{lem:edge_remove_approx}
Suppose $X \in \Omega$, $e = (u,v)$, and $p^{*} \in [0,1]$ are such that $\sup_{u\in [n]}|p_u(X)-p^{*}| \leq \epsilon$. 
% Fix a graph $G$ and let $d_{ij}=d_{ij}(G)$ as defined in \eqref{e:pairDeg}. Then 
For any fixed graph $G$ and for some constant $C_G$ depending only on $G$, we have
$$\Big|\delN e G X - \frac{2}{n^{2}}\sum_{(i,j)\in E(G)}\B(\frac{p_{uv}(X)}{(\ps)^2}\B)^{d_{ij}}(p^{*})^{|E(G)| -1}\Big| \leq \frac{C_{G}}{n^2}\left(\epsilon + \delta_{\square}(\tilde{X},p^{*}) + n^{-1}\right)\,. $$
\end{lemma}

\begin{corollary}\label{cor:itemsUnif}
Let $C(L)=\max_{G\in  \mathbb{G}_L} C_G$ for $C_G$ in Lemma~\ref{lem:edge_remove_approx} and $h(\eps,L)$ be as in \eqref{eq:delt_conc}. Let $a(\eps,L) = h(\eps,L)/3C(L)$.
Then $ X \in \Gam {\ps} \epsilon$ holds if
\begin{enumerate}
    \item $X$ is $a(\eps,L)$-close to the constant graphon $p^{*}$ in the cut-metric,
    \item $p_u(X)$ is uniformly close to $p^{*}$ for every vertex $u$, i.e., $\sup_u |p_u(X) - \ps|\leq a(\eps,L)$, and
    \item $p_{uv}(X)$ is close to $(p^{*})^2$ uniformly for every pair of vertices $u \neq v$, i.e., 
    \begin{equation}\label{eq:doub_req}
\sup_{\substack{ u\neq v}}|p_{uv}(X)-(p^{*})^2| \leq \ps a(\eps,L)/3L\,.
\end{equation} It follows that Theorem~\ref{thm:final_unif_conc} is proved if these conditions are each shown to hold for $X\sim \mu(\cdot|\Ball \ps \eta )$
with probability at least $1- C(\beta,\epsilon,\eta,L)\exp(-\Omega_{\beta,\epsilon,\eta,L}(n))$
\end{enumerate}
\end{corollary}
Item 1 holds by Theorem~\ref{thm:large_deviations}. 
We address Items 2 and 3 below.
% Item 2 will  holds with high probability by Corollary~\ref{cor:degree_conc}. Therefore, the proof of Theorem~\ref{thm:final_unif_conc} reduces to showing that Item 3 holds with high probability under the measure $\mu\b(\cdot\b|\Ball \ps \eta\b)$. 

\subsection{Uniform Control of Degrees}

Section~\ref{sec:cavity} develops the cavity method for the $\ergm$ and demonstrates the following uniform control on vertex degrees.

% \begin{customthm}{8}\label{eight}
% Every theorem must be numbered by hand.
% \end{customthm}

 \begin{customthm}{\ref*{cor:degree_conc}}%\label{cor:degree_conc}
 Suppose $p^{*} \in U_{\beta}$ and let
$\epsilon > 0$ be an arbitrary fixed constant. Then, we can take $0<\eta < c(\beta,\epsilon) $ and $ n > n_0(\beta,\epsilon,\eta)$ such that
$$\mu\Big(\sup_{u \in [n]}|p_u(X) - p^{*}| \leq \epsilon\,\Big| \Ball \ps \eta\Big) \geq 1-\exp(-\Omega_{\beta,\epsilon}(n))\,.$$
 \end{customthm}

 We 
%  replace $\epsilon$ in Lemma~\ref{lem:edge_remove_approx} with $h(\epsilon,L)$ as used in Equation~\eqref{eq:delt_conc}, and
 take 
 $\eta< c(\beta, a(\eps,L)$
%  $\eta \leq \min\{h(\eps,L),c(\beta,h(\eps,L))\}$ 
 in Corollary~\ref{cor:degree_conc}, implying that if we sample $X \sim \mu\b(\cdot\b|\Ball \ps \eta\b)$ then Item 2 of Corollary~\ref{cor:itemsUnif}
 holds with probability $ 1- C(\beta,\epsilon,\eta,L)\exp(-\Omega_{\beta,\epsilon,\eta,L}(n))$. 
%  If $1/n \leq h(\epsilon,L)$.
%  we conclude from Lemma~\ref{lem:edge_remove_approx} that Equation~\eqref{eq:delt_conc} holds with probability at least $1- C(\beta,\epsilon,\eta,L)\exp(-\Omega_{\beta,\epsilon,\eta,L}(n))$, whenever
% \begin{equation}\label{eq:doub_req}
% \sup_{\substack{u,v \\ u\neq v}}|p_{uv}(X)-(p^{*})^2| \leq c_{\beta,L}h(\epsilon,L)
% \end{equation} with probability at least $1- C(\beta,\epsilon,\eta,L)\exp(-\Omega_{\beta,\epsilon,\eta,L}(n))$ for some small enough $c(\beta,L)$. Thus, if~\eqref{eq:doub_req} holds with high probability then we have proved Theorem~\ref{thm:final_unif_conc}. 

% We do this below. 

% Recall that in Lemma~\ref{lem:almost_degree_conc}, we showed that whenever $X$ is close to the constant graphon $p^{*}$, then \emph{most} of the degrees concentrate close to $p^{*}$. 
% In order to apply Lemma~\ref{lem:edge_remove_approx}, we also require control of common neighbor $p_{uv}(X)$ counts. 
It remains to show Item 3 of Corollary~\ref{cor:itemsUnif}.
As a step towards this, it turns out that if $p_{u}(X)$ is close to $p^{*}$ uniformly for \emph{every} vertex $u$ (as stated in Corollary~\ref{cor:degree_conc}) and $X$ is close to the constant graphon $p^{*}$, then \emph{most} of the common neighbor counts $p_{uv}(X)$ are close to $(p^{*})^2$. This follows from the definition of cut-metric. 
% That is, if the degrees are uniformly close to $p^{*}$ but \emph{most} $p_{uv}(X)$ are not close to $(p^{*})^2$, then we can find a cut such $\delta_\square(\tilde{X},p^{*})$ is large which contradicts the assumption that $\delta_\square(\tilde{X},p^{*})$ is small. 
% This will be useful in Section~\ref{sec:unif_conc}, where we will prove Theorem~\ref{thm:final_unif_conc}, which is the important component behind the proof of Theorem~\ref{thm:main_result}.

%\gb{Explain why this next thing is true}

\begin{customthm}{\ref*{cor:almost_dub_conc}}
Suppose that $p^{*}\in U_{\beta}$. Given arbitrary $\epsilon, \delta > 0$, suppose $\eta < c(\beta,\epsilon,\delta)$ and $n > n_0(\delta,\epsilon,\beta,\eta)$. Then, for every $ u\in [n]$, there exists a random set $S_u \subseteq [n]\setminus \{u\}$ such that $|S_u|\leq \delta n$ and
$$
\mu\Big(\sup_{u\in [n]}\sup_{v\in S^{\complement}_u}|p_{uv}(X)-(p^{*})^2| \leq \epsilon\,\Big|\Ball \ps\eta\Big) \geq 1-\exp(-\Omega_{\beta,\epsilon}(n))\,.
$$
\end{customthm}
We refer to Appendix~\ref{sec:cor_proof_2} for the proof. 

Next, to establish Equation~\eqref{eq:doub_req}, we boost control of $p_{uv}(X)$ from most pairs $u,v$ to all pairs.

\subsection{Uniform Control of Common Neighbors}
Let $X^{\prime}$ be obtained from $X \sim \mu\b(\cdot\b|\Ball \ps \eta\b)$ via one step of the Glauber dynamics with respect to $\mu\b(\cdot\b|\Ball \ps \eta\b)$, so that we also have $X' \sim \mu\b(\cdot\b|\Ball \ps \eta\b)$. As shown in Lemma~\ref{lem:outer_prob}, 
with high probability $X$ is $\eta/2$ away from the boundary of $\Ball \ps \eta$ and the expected Glauber update for $p_{uv}$ with respect to $\mu\b(\cdot\b|\Ball \ps \eta\b)$ is the same as with respect to $\mu(\cdot)$, which is
\begin{align*}
&\mathbb{E}\left[p_{uv}(X^{\prime})\b|X\right] \\&\quad = \B(1-\frac{2}{N}\B)p_{uv}(X) + \frac{1}{N}\B(n\inv\sum_{w \in [n]\setminus\{u,v\}}\phi_{uw}(X_{\sim uw})X_{vw} + n\inv\sum_{w \in [n]\setminus\{u,v\}}\phi_{vw}(X_{\sim vw})X_{uw}\B) \,.
\end{align*} 
The next lemma shows that under 
the conditions shown in Corollaries~\ref{cor:degree_conc} and~\ref{cor:almost_dub_conc} to hold with high probability for $X \sim \mu\b(\cdot\b|\Ball \ps \eta\b)$,
each of the normalized sums in the last displayed equation is close to $(p^{*})^2$
% the quantity $$n\inv\sum_{w \in [n]\setminus\{u,v\}}\phi_{uw}(X_{\sim uw})X_{vw} \approx (p^{*})^2\,,$$ 
uniformly for every $u\neq v$. The Glauber dynamics with respect to $\mu\b(\cdot\b|\Ball \ps \eta\b)$, therefore, tries to regress every $p_{uv}(X)$ close to $(p^{*})^2$.

\begin{lemma} \label{lem:fixed_pt_crude}
Suppose $X \in \Omega$, $p^{*} \in [0,1]$ are such that the following conditions hold:
\begin{enumerate}
    \item $\sup_{u\in [n]}|p_u(X)-p^{*}| \leq \epsilon$, and
    \item For every $u \in [n]$, there exist sets $S_{u} \subseteq [n]$ such that $|S_u| \leq \delta n$ and
$$\sup_{u\in [n]}\sup_{w\in S^{\complement}_u}|p_{uw}(X)-(p^{*})^2| \leq \epsilon\,.$$
\end{enumerate} 
Then, for every $u,v \in [n]$, we have
$$\Big|n\inv
% \sum_{\substack{w \in [n]\\ w \notin \{u,v\}}}
\sum_{w \in [n]\setminus\{u,v\}}
\phi_{uw}(X_{\sim uw})X_{vw} -(p^*)^2\Big| \leq C_{\beta}\left(\epsilon + \delta +\delta_{\square}(X,p^{*}) + n\inv\right)\,.$$
\end{lemma}

The lemma is proved in the next subsection.

Intuitively, this suggests that the stationary distribution of the Glauber dynamics, $\mu\b(\cdot\b|\Ball \ps \eta\b)$, should be such that $p_{uv}(X) \approx (p^{*})^2$ for every $u \neq v$ with high probability. The next 
% Lemma~\ref{lem:fixed_pt_conc} along with Lemma~\ref{lem:fixed_pt_crude} 
lemma formalizes this sentiment using Stein's method for concentration developed in \cite[Theorem 1.5]{chatterjee2007stein}.

\begin{lemma}\label{lem:fixed_pt_conc}
For any $u,v \in [n]$ such that $u \neq v$,
define $$g_{uv}(X) := p_{uv}(X) - \frac{1}{2n}
% \sum_{\substack{w \in [n]\\ w \notin \{u,v\}}}
\sum_{w \in [n]\setminus\{u,v\}}
\B(\phi_{uw}(X_{\sim uw})X_{vw}+\phi_{vw}(X_{\sim vw})X_{uw}\B)\,.$$

Then, for any $\gamma > 0$ which is independent of $n$, we have 
\begin{equation}
    \mu\Big( \{|g_{uv}(X)| > \gamma\}\cup (\Ball\ps{\eta/2})^{\complement}  \,\Big| \Ball \ps \eta\Big) \leq C(\eta,\beta)\exp(-\Omega_{\gamma,\eta,\beta}(n))\,.
\end{equation}
\end{lemma}

\begin{corollary}
Item 3 of Corollary~\ref{cor:itemsUnif} holds with the desired probability. 
\end{corollary}
\begin{proof}
Combining the last two lemmas with Corollary~\ref{cor:degree_conc} proves \eqref{eq:doub_req}.
\end{proof}

\subsection{Proof of Lemma~\ref{lem:fixed_pt_crude}}

First, note that by definition
$$ \phi_{uw}(X_{\sim uw}) = \frac{\exp(\sum_{i=0}^{K}n^2\beta_i \delN {uw} i X)}{1+\exp(\sum_{i=0}^{K}n^2\beta_i \delN {uw} i X)}\,.$$ Now, suppose $w \in S^{\complement}_u$. Then, using Lemma~\ref{lem:edge_remove_approx}, we conclude that $$\sum_{i=0}^{K}n^{2}\beta_i \delN {uw} i X = \sum_{i=0}^{K}2|E_i|\beta_i(p^{*})^{|E_i|-1} \pm O_{\beta}\left(\epsilon + \delta_{\square}(X,p^{*})+{n\inv}\right)\,.$$ Using the fact that the function $x \to \frac{e^x}{1+e^x}$ is $1$-Lipschitz, we have for $w \in S^{\complement}_u$ that
$$\phi_{uw}(X_{\sim uw}) = \frac{\exp(\sum_{i=0}^{K}2|E_i|\beta_i (p^{*})^{|E_i|-1})}{1+\exp(\sum_{i=0}^{K}2|E_i|\beta_i (p^{*})^{|E_i|-1})} \pm O_{\beta}\left(\epsilon + \delta_{\square}(X,p^{*})+n\inv\right)\,.$$
Now, the fact that $p^{*} \in U_\beta$ implies that $L_{\beta}^{\prime}(p^{*}) = 0$. It can be easily checked that this implies $p^{*} = \phi_\beta(\ps) $. It follows that whenever $w \in S^{\complement}_u$,
$$\phi_{uw}(X_{\sim uw}) = p^{*} \pm O_{\beta}\left(\epsilon + \delta_{\square}(X,p^{*})+n\inv\right)\,.$$

An application of the triangle inequality now shows that
$$\B|n\inv\sum_{\substack{w \in [n]\setminus \{u,v\}}}\phi_{uw}(X_{\sim uw})X_{vw} - p^{*}p_v(X)\B| \leq \frac{|S_u|}{n} + O_{\beta}\left(\epsilon + \delta_{\square}(X,p^{*})+n\inv\right)\,.
$$
The assumption that $|S_u| \leq \delta n$ and the fact that $|p_v(X) - p^{*}| \leq \epsilon $ imply the
result.
\qed

%% file: technical.tex
\section{Some Graphon Estimates}

In this section, we state several technical results that will be needed later. Results from other works are stated without proof and otherwise the proofs are given in the appendix. 

% \gb{I think we may want to move 6.1 up to the Graphon section or large deviations section}

% First, we state the following results obtained by combining Theorem 3.2, Theorem 4.1, and Theorem 4.2 in \cite{chatterjee2013estimating}. Let $U_{\beta}$ be the finite set of global maximizers of $L_{\beta}$ as defined in Section~\ref{sec:definition}. Define the set $\tilde{U}_{\beta} := \{p \mathbf{1} : p \in U_{\beta} \} \subset \Tilde{\mathcal{W}}$.

% \begin{theorem}
% \label{thm:large_deviations}
% Let $X_{n} \sim \ergm(n,\beta)$ and $\tilde{X}_n$ be its corresponding graphon. For every fixed $\eta > 0$, there are constants $C(\eta),c(\eta) > 0$ such that
% $$\mathbb{P}(\delta_\square(\tilde{X}_n,\tilde{U}_{\beta})> \eta) \leq C(\eta)\exp(-c(\eta)n^2)\,.$$
% \end{theorem}

Even though Theorem 3.2 in \cite{chatterjee2013estimating} considers the probability over the entire space of graphons, we can easily adapt its proof to show the following lemma which considers only the neighborhood $\Ball \ps {\eta+\delta}$ for $\ps \in U_{\beta}$.

\begin{lemma}
\label{lem:outer_prob}
Suppose $p^{*} \in U_{\beta}$. Then, there exists a constant $c_{\beta} > 0$ such that whenever $\eta,\delta \in (0,c_{\beta})$ are fixed constants independent of $n$, we have
$$\frac{\mu(\eta \leq \delta_\square(X,p^{*}) \leq \eta+\delta)}{\mu(\delta_\square(X,p^{*}) \leq \eta)} \leq C(\eta,\delta)\exp(-c(\eta,\delta)n^2)\,.$$
\end{lemma}

% Next, we state the classic \emph{inverse counting lemma} from the graphon literature. We refer to \cite[Lemma 10.32]{lovasz2012large} for the more general result.

% \begin{lemma}
% \label{lem:inv_count}
% Given any $\epsilon > 0$, there exists $k,\eta >0$ such that for any graphons $h_1,h_2$ if \gb{made some edits, check}
% $$\sup_{G:  |V(G)|\leq k}|t(G,h_1)-t(G,h_2)| \leq \eta\,,$$
% then
% $\delta_{\square}(h_1,h_2) \leq \epsilon$.
% \end{lemma}
% \dn{might not need the lemma above, check again!}

We state three technical lemmas below, whose proofs appear in Appendix~\ref{sec:deferred_proofs}.

\begin{lemma}
\label{lem:almost_degree_conc}
For any graph $X$ over $n$ vertices, $p \in [0,1]$, and $\delta > 0$, there exists a set $S\subseteq [n]$ such that $|S| \leq \delta n$ and 
\begin{equation}\label{eq:edge_cut_metric}
    \sup_{u\in S^{\complement}} |p_u(X) - p| \leq  \frac{2\delta_{\square}(\tilde{X},p\mathbf{1}) }{\delta}\,.
\end{equation}
It follows as an application that for $X\sim \ergm(n,\beta)$ and any given fixed constants $\delta,\eta > 0$, with probability at-least $1-C(\eta)\exp(-c(\eta)n^2)$, there exists a (random) set $S \subseteq [n]$ such that $|S| \leq \delta n$ and
$$\inf_{p^{*}\in U_{\beta}}\sup_{u\in S^{\complement}} |p_u(X) - p^{*}| \leq  \frac{2\eta }{\delta}\,.$$
\end{lemma}

\begin{lemma} \label{lem:dom_bound_delta}
% Given two graphs $X,Y \in \Omega$, we denote by $X \preceq Y$ the relation $X_e \leq Y_e$ for every $e \in {{[n]}\choose{2}}$. 
If $Y \preceq X \preceq Z$, then $\delta_{\square}(\tilde{X},p\mathbf{1}) \leq \max\big(\delta_{\square}(\tilde{Y},p\mathbf{1}), \delta_{\square}(\tilde{Z},p\mathbf{1})\big) $.
\end{lemma}

% Consider any fixed graph $G = (V,E)$. For any vertex $i$, let $d_i$ denote the degree of the vertex $i$. 

% \begin{lemma}\label{lem:young_ineq}
% Suppose $p,q \in [0,1]$. Then
% $$\sum_{i=1}^{|V|}q^{d_i}p^{|E|-d_i} \leq 2q^{|E|} + (|V|-2)p^{|E|}\,.$$
% This is an equality when the graph is a single edge connecting two vertices (which we denoted by $G_0$). Now suppose that the graph $G$ is connected and $ G \neq G_0$. Then we have 
% \begin{align*}
%     &2q^{|E|} + (|V|-2)p^{|E|} - C|p-q|
%     \\&\qquad\qquad\leq \sum_{i=1}^{|V|}q^{d_i}p^{E-d_i} \leq 2q^{|E|} + (|V|-2)p^{|E|} - \zeta(|p-q|)
% \end{align*}\,.
% Here $\zeta: [0,1] \to \mathbb{R}^{+}$ is a continuous function depending only on $G$ such that $\zeta(x) > 0$ whenever $x \neq 0$ and $C$ is a constant depending only on $(d_i)_{i \in V}$ and $|E|$. 
% \end{lemma}

\begin{lemma}\label{lem:young_ineq}
Consider any fixed graph $G = (V,E)$. For any vertex $i$, let $d_i$ denote the degree of the vertex $i$. 
Suppose $p,q \in [0,1]$. 
If $G=G_0$, the graph consisting of a single edge, then
\begin{equation}\label{eq:holder_1}
\sum_{i=1}^{|V|}q^{d_i}p^{|E|-d_i} = 2q^{|E|} + (|V|-2)p^{|E|}\,.\end{equation}
If $G$ is connected and $ G \neq G_0$, then there is a constant $C$ depending only on $(d_i)_{i \in V}$ and $|E|$ such that
\begin{align*}
    &2q^{|E|} + (|V|-2)p^{|E|} - C|p-q|
    \\&\qquad\qquad\leq \sum_{i=1}^{|V|}q^{d_i}p^{E-d_i} \leq 2q^{|E|} + (|V|-2)p^{|E|} - \zeta(|p-q|)\,,
\end{align*}
where $\zeta: [0,1] \to \mathbb{R}^{+}$ is a continuous function depending only on $G$ such that $\zeta(x) > 0$ for $x \neq 0$. 
\end{lemma}

%% file: cavity.tex
\section{The Cavity Method}\label{sec:cavity}

In this section we address the degrees and show that every vertex $u$ has nearly the same degree $p_u(X) \approx p^{*} \in U_{\beta}$ with high probability for $X \sim \mu(\spacedot|\Ball \ps \eta)$. 
While the cut-metric based convergence does not allow us to control all the degrees, it is nevertheless possible to conclude that a large portion of the vertices have degree $p_u(X)\approx \ps$.
We boost this to a uniform statement, in Theorem~\ref{thm:cavity_conc} and Corollary~\ref{cor:degree_conc}, via the \emph{cavity method}: most of the vertices and the corresponding edges are conditioned on being close to the constant graphon $p^{*}\mathbf{1}$, which generates the \emph{mean field} with which the remaining cavity vertices interact. We can then reason about the behavior of the cavity vertices.
% Considering the distribution of the cavity after conditioning on the mean field being close to the constant graphon $p^{*}\mathbf{1}$ 
% This approach yields the main results of this section, Theorem~\ref{thm:cavity_conc} and Corollary~\ref{cor:degree_conc}, which are the desired uniform vertex concentration.
% Theorem~\ref{thm:cavity_conc} 
% shows that the probability of the degree of any cavity vertex being far away from $p^{*}$ has exponentially small probability. Corollary~\ref{cor:degree_conc} then shows that $p_u(X)$ indeed concentrates near $p^{*}$ uniformly for every $u \in [n]$. 

We start by adapting several of the graphon definitions to incorporate a cavity.
% Before delving into the cavity method, we will first carry out some useful calculations regarding graphons. 

% Their proofs are straightforward and can be found in Appendix~\ref{sec:deferred_proofs}. 
\subsection{Restricted Homomorphism Densities and Restricted Cut Metric}

Recall from Section~\ref{s:graphon} the function representative $f^X(x, y)$ of a graph $X$ over $n$ vertices. We will need the homomorphism density of a graph forced to contain a particular vertex $u$ of $X$. 
To that end, for every $u \in [n]$, define the event
$$A^{k,i}_u := \B\{x \in [0,1]^k: x_i \in  \B
[\frac{u-1}{n},\frac{u}{n}\B) \B\}$$
and
\begin{equation}\label{e:Aku}
    A^{k}_u = \B\{x \in [0,1]^k: x_i \in  \B[\frac{u-1}{n},\frac{u}{n}\B) \text{ for some } i \in [k] \B\}=\bigcup_{i\in [k]}A^{k,i}_u\,.
\end{equation}

For the sake of clarity in the results below, given a fixed graph $G$, we take its vertex set $V(G) = [k]$ during calculations. 

\begin{definition}[Homomorphism density w.r.t. a vertex]
Define the homomorphism density of $G$ in $X$ with respect to vertex $u$ (which counts only homomorphisms which include the vertex $u$) as
\begin{equation}
    N_G(X;u) := \int_{[0,1]^{k}} \mathbbm{1}(A^{k}_u)\prod_{(i,j)\in E(G)}f^{X}(x_i,x_j) \prod_{i=1}^{k}dx_i\,.
\end{equation}
For ease of computation, we also introduce the quantity
\begin{equation}
    N^{0}_G(X;u) := \sum_{l=1}^{k}\int_{[0,1]^{k}} \mathbbm{1}(A^{k,l}_u)\prod_{(i,j)\in E(G)}f^{X}(x_i,x_j) \prod_{i=1}^{k}dx_i\,.
\end{equation}
\end{definition}

The next lemma follows from elementary 
arguments; see Appendix~\ref{sec:deferred_proofs}.
\begin{lemma}\label{lem:inner_count}
$0 \leq N_G^{0}(X,u) - N_G(X,u) \leq {|V(G)|^3}/{n^2}$.
\end{lemma}
We will now show that whenever the graphon corresponding to $X$ is close to a constant graphon $p\mathbf{1}$, then $N_G(X,u)$ can be approximated as a polynomial of $p_u(X)$ and $p$. This will allow us to control the fine-grained structure of $X$ in terms of the counts $N_G(X,u)$ just based on nearness to a constant graphon and the normalized degrees of the vertices of $X$. The proof, given in Appendix~\ref{sec:deferred_proofs}, follows from a slight modification of the proof technique of 
\cite[Lemma 4.4]{borgs2006graph}, which establishes the continuity of homomorphism densities with respect to the graphon metric via a repeated application of the triangle inequality.

\begin{lemma}\label{lem:single_vertex_decomp}
 Suppose $D_G = (d_1,\dots,d_{|V(G)|})$ is the degree sequence of the fixed graph $G$ considered above. For any graph $X$ with vertex set $[n]$ and $u \in [n]$, we have
 $$\B|N_G(X;u) -n\inv\sum_{d\in D_G}p_u(X)^{d}p^{|E(G)|-d}\B| \leq |V(G)||E(G)|\frac{\delta_\square(\tilde{X},p)}{n}  +\frac{|V(G)|^3}{n^2}\,.$$
\end{lemma}

% \subsection{Cavity Decomposition}

 Let $S\subset [n]$ be the ``cavity set". Define $$A_S^{k} := \cup_{u\in S}A_u^{k} \quad \text{and}\quad A_S^{k,i} := \cup_{u\in S} A_u^{k,i}\,.$$ 
 We have that 
\begin{align}
    N_G(X) &= N_G(X;S)+\int_{[0,1]^k} \mathbbm{1}\big\{(A_S^{k})^\complement
    % A_S^{k,\complement}\big) 
    \big\}
    \prod_{(i,j)\in E(G)} f^{X}(x_i,x_j)\prod_{i=1}^{k}dx_i \nonumber \,,
    % \\ &\qquad \quad +  \int_{[0,1]^k} \mathbbm{1}\big\{A_S^{k}\big\} \prod_{(i,j)\in E(G)} f^{X}(x_i,x_j)\prod_{i=1}^{k}dx_i \,.
\end{align}
% We will denote the second term by
where the subgraph count $N_G(X;S)$ restricts to subgraphs containing a cavity vertex,
$$N_G(X;S) := \int_{[0,1]^k} \mathbbm{1}\big\{A_S^{k}\big\} \prod_{(i,j)\in E(G)} f^{X}(x_i,x_j)\prod_{i=1}^{k}dx_i\,.$$

The proofs of the next two lemmas are deferred to Appendix~\ref{sec:deferred_proofs}.
\begin{lemma}\label{lem:inner_event_S}
Let $X$ be any simple graph with vertex set $[n]$ and let $S \subset [n]$ be arbitrary. Then
$$\Big|N_G(X;S) - \sum_{u\in S}N_G(X;u)\Big| \leq \frac{k^3|S|^2}{n^2}\,.$$
\end{lemma}

 We now define the graphon metric restricted to $S^{\complement}$.

 \begin{definition}[Restricted Graphon Metric]\label{def:res_met}
  Let $p^{*}\in [0,1]$. Define \begin{equation}
    f^{X,S,p^{*}}(x_1,x_2) = \begin{cases}
    f^{X}(x_1,x_2) \text{ if } \lceil nx_1\rceil,\lceil n x_2\rceil \in S^{\complement} \\
    p^{*} \text{ otherwise}\,,
    \end{cases}
\end{equation}
and let $\tilde{X}^{S,p^{*}}$ be the graphon corresponding to $f^{X,S,p^{*}}$. We define the restricted cut metric to be $$\delta_{\square}^{S,p^{*}}(\tilde{X},p^{*}) := \delta_{\square}(\tilde{X}^{S,p^{*}},p^{*})\,.$$ 
 \end{definition}
 
The restricted graphon distance to $p^*$ can be approximated in terms of the unrestricted distance:

\begin{lemma}\label{lem:restricted_metric}
For any $S \subseteq [n]$, we have
$$ \delta_{\square}(\tilde{X},p^{*})- \frac{|S|(2n-|S|)}{n^2}\leq \delta_{\square}^{S,p^{*}}(\tilde{X},p^{*}) \leq \delta_{\square}(\tilde{X},p^{*})\,.$$
\end{lemma}

We are now ready to establish the cavity decomposition of the Hamiltonian.

\subsection{Cavity Decomposition of the Hamiltonian}
%\gb{Should we use $s$?}
It will be convenient to let $r := |S|$. Given a simple graph $X$ over $n$ vertices, we define $\bar{p}_u(X)$ to be the number of edges from vertex $u \in [n]$ to the set $S^{\complement}$, normalized by $n$:
$$\bar{p}_u(X) = \frac1{n}{\sum_{v \in S^{\complement}}X_{uv}}\,.$$
Here and throughout we hide the dependence on $S$ to streamline the notation. Additionally, whenever it is clear, we will denote $\bar{p}_u(X)$ by $\bar{p}_u$.
Note that $|\bar{p}_u(X)-p_u(X)| \leq {|S|}/{n}$. 

% We now present the cavity decomposition of the Hamiltonian. 
Denote the portion of the Hamiltonian associated to the cavity by
$$
\Hcav(X;S) := \sum_{i=0}^{K}n^2\beta_i N_i(X;S)\,,
$$ 
which is the same as $\mathcal{H}_\beta$ except that the homomorphism densities are restricted to have at least one vertex in the set $S$. Denote the rest of the Hamiltonian by
$$
\Hmean(X;S) = \mathcal{H}_{\beta}(X)- \Hcav(X;S)\,.
$$
We next bound the difference between $\mathcal{H}_\beta(X)$ and $\Hmean(X;S)$.

\begin{lemma}[Cavity Decomposition]
\label{lem:cavity_decomp}
Assume that $\beta_i >0$ for some $ i = 1,\dots, K$. We have the following upper and lower bounds.
\begin{enumerate}
    \item Upper Bound:
    \begin{align}\mathcal{H}_{\beta}(X) &\leq \Hmean(X;S) + O_{\beta}(nr\delta^{S,p^{*}}_\square(\tilde{X},p^{*}) + r^2) + \sum_{i=0}^{K}\beta_i n r |V_i|(p^{*})^{|E_i|} 
 \nonumber \\ &\qquad+2n\sum_{u\in S}\Big[ L_{\beta}(\bar{p}_u)+I(\bar{p}_u)-L_{\beta}(p^{*})-I(p^{*})- \zeta_{\beta}(|\bar{p}_u-p^{*}|)\Big] \end{align}
 \item Lower Bound:
  \begin{align}\mathcal{H}_{\beta}(X) &\geq \Hmean(X;S) - O_{\beta}(nr\delta^{S,p^{*}}_\square(\tilde{X},p^{*}) + r^2) + \sum_{i=0}^{K}\beta_i n r |V_i|(p^{*})^{|E_i|} 
 \nonumber \\ &\qquad+2n\sum_{u\in S}\Big[ L_{\beta}(\bar{p}_u)+I(\bar{p}_u)-L_{\beta}(p^{*})-I(p^{*})- C_{\beta}|\bar{p}_u-p^{*}|\Big] \end{align}
\end{enumerate}
Here $\zeta_{\beta}$ is a function with the same properties of the function $\zeta$ in Lemma~\ref{lem:young_ineq}
\end{lemma}
\begin{proof}
Let $D_i$ be the tuple of the degrees of vertices in $G_i$ (as in Lemma~\ref{lem:single_vertex_decomp}). In the equations below, we will take $\sum_{d_i\in D_i}$ to mean summation over all elements of the tuple. We have
\begin{align}
    \mathcal{H}_{\beta}(X) &= \Hmean(X;S)+ \Hcav(X;S)\nonumber \\ &= \Hmean(X;S) + \sum_{u\in S}\sum_{i=0}^{K}\beta_i n^2N_i(X;u) \pm O_{\beta}(r^2) \nonumber \\
&= \Hmean(X;S) + \sum_{u\in S}\sum_{i=0}^{K}\sum_{d_i\in D_i}\beta_i n p_u^{d_i}(p^{*})^{E_i-d_i} \pm  O_{\beta}(nr\delta_\square(\tilde{X},p^{*}) + r^2) \nonumber \\
&= \Hmean(X;S) + \sum_{u\in S}\sum_{i=0}^{K}\sum_{d_i\in D_i}\beta_i n \bar{p}_u^{d_i}(p^{*})^{E_i-d_i} \pm  O_{\beta}(nr\delta_\square(\tilde{X},p^{*}) + r^2)\,. \label{eq:cavity_simplify_1}
\end{align}
The first step is by the definition of the cavity Hamiltonian. The second step uses Lemma~\ref{lem:inner_event_S} to approximate $N(G;S)$ and the third step uses Lemma~\ref{lem:single_vertex_decomp} to approximate $N_i(X;u)$. In the fourth step, we have used the fact that $|\bar{p}_u - p_u|\leq {r}/{n}$. 

We now apply Lemma~\ref{lem:young_ineq} to the second term of Equation~\eqref{eq:cavity_simplify_1}, yielding
\begin{align}
&\sum_{i=0}^{K}\sum_{d_i\in D_i} \beta_i\bar{p}_u^{d_i}(p^{*})^{E_i-d_i} \leq \sum_{i=0}^{K}\beta_i\Big[ 2\bar{p}_u^{|E_i|} + (|V_i|-2)(p^{*})^{|E_i|}\Big] - \zeta_{\beta}(|\bar{p}_u - p^{*}|) \nonumber \\ 
&= \sum_{i=0}^{K}\beta_i  |V_i|(p^{*})^{|E_i|} 
+2\Big[ L_{\beta}(\bar{p}_u)+I(\bar{p}_u)-L_{\beta}(p^{*})-I(p^{*})- \zeta_{\beta}(|\bar{p}_u-p^{*}|)\Big]\,. \label{eq:apply_young}
\end{align}
In the first step we have used Lemma~\ref{lem:young_ineq} and the fact that for $G_0$ the inequality is an equality which allows for all $\beta_0 \in \mathbb{R}$. For $i > 0$, notice that $\beta_i \geq 0$ and the inequality goes in the right direction. The function $\zeta_{\beta}$ is as defined in the statement of Lemma~\ref{lem:young_ineq} and exists since $\beta_i > 0$ for some $i \in [K]$.  To see this, observe that the Equation~\ref{eq:holder_1} is an equality when $G=G_0$. Therefore, in order to establish the strict inequality involving $\zeta$ as shown in Lemma~\ref{lem:young_ineq}, we need at-least one of the $\beta_i > 0$. 

The upper bound in the lemma statement follows by combining Equations~\eqref{eq:cavity_simplify_1} and~\eqref{eq:apply_young} along with Lemma~\ref{lem:restricted_metric} to show that $$nr\delta_\square(\tilde{X},p^{*}) \leq nr \delta^{S,p^{*}}_\square(\tilde{X},p^{*}) + 2r^2 \,.$$ The lower bound on the Hamiltonian follows from a similar argument by replacing the upper bound in Lemma~\ref{lem:young_ineq} with the lower bound. 
\end{proof}

\subsection{Controlling Degrees of Cavity Vertices}

Given a sequence $\bar{q}_u  \in \{0,1/n,2/n,\dots,1\}$ for $u \in S$, we define the tuple $\mathbf{q}_S = (\bar{q}_u)_{u\in S}$. Given arbitrary and fixed $\eta > 0$, $p^{*} \in [0,1]$, we define the events
\begin{equation}
A(S,\mathbf{q}_S,p^{*},\eta) = \{X: \bar{p}_u(X) = \bar{q}_u \text{ for } u \in S\}\cap\{ \delta^{S,p^{*}}_\square(\tilde{X},p^{*}) \leq \eta \}
\end{equation}
and
\begin{equation}
\label{eq:B_def}
B(S,p^{*},\eta) := \{\delta^{S,p^{*}}_\square(\tilde{X},p^{*}) \leq \eta\}\,.
\end{equation}
Note that by definition, $B(\emptyset,p^{*},\eta) = \Ball \ps \eta$. We want to show that whenever $p^{*} \in U_{\beta}$, if $\mathbf{q}_S$ is not close to $p^{*}$, then the event $A(S,\mathbf{q}_S,p^{*},\eta)$ has exponentially small probability compared to the event $B(S,p^{*},\eta)$, whenever $\eta$ and $S$ are small enough. 

We now note that $\Hmean(X;S)$ and $\delta^{S,p^{*}}_\square(\tilde{X},p^{*})$ depend only on $X_{uv}$ for $u,v \in S^{\complement}$. Therefore, whenever $|S|$ is small, we will think of $\Hmean(X;S)$ as the mean field which controls the behavior of the cavity, i.e., the edges emanating from the vertices in $S$. Now, fixing $X$ such that $\delta^{S,p^{*}}_\square(\tilde{X},p^{*}) \leq \eta$, we look at the joint law of $(X_{uv})$ such that at least one of $u$ or $v$ is in the set $S$. By $\Xmean(S)$ we denote the coordinates $(X_{uv})_{u,v \in S^{\complement}}$. We denote the rest of the coordinates by $\Xcav(S)$. Therefore, we want to understand the conditional law $\Xcav(S)|\Xmean(S)$ under the measure $\mu$. We first record the following combinatorial lemma, whose proof can be found in Appendix~\ref{sec:deferred_proofs}.

\begin{lemma}\label{lem:count_approx} Suppose ${r}/{n}\leq {1}/{2}$
and $\Xmean(S)$ is fixed. Let the count of $\Xcav (S)$ such that $\bar{p}_u(X) = q_u$ for $u \in S$ be denoted by $H_{\mathsf{cav}}(\mathbf{q}_{S})$. $H_{\mathsf{cav}}(\mathbf{q}_{S})$ satisfies
$$\exp\left( - r^2\left[4+2\log(\tfrac{n}{r})\right]-\frac{r}{2}\log(2n) \right)\leq H_{\mathsf{cav}}(\mathbf{q}_S)\exp\Big(2n\sum_{u\in S}I(q_u)\Big)\leq 1
\,.$$ 
\end{lemma}
%\gb{Should we use diff notation for this given that most of the time $N$ is used for something else?}

Below we present the main result of this section. 
\begin{theorem}\label{thm:cavity_conc}
Let $p^{*} \in U_{\beta}$ be such that ${r}/{n} < p^{*} < 1-{r}/{n}$ and $r < {n}/{2}$. Given any $\mathbf{q}_S = (\bar{q}_u)_{u\in S}$ as defined above for $n$ sufficiently large as a function of $\beta,\eta$, we have
\begin{align*}\frac{\mu(A(S,\mathbf{q}_S,p^{*},\eta))}{\mu(B(S,p^{*},\eta))} 
 \leq \exp\bigg(2n\sum_{u\in S} \Big(L_{\beta}(\bar{q}_u)-L_{\beta}(p^{*})-\zeta_{\beta}(|\bar{q}_u - p^{*}|)\Big)+ \text{lower order}\bigg)\,,\end{align*} 
where 
$\text{lower order} = O_{\beta}(nr\eta + r^{2}\log({n}/{r})+ r\log n)$.
\end{theorem}

\begin{proof}
Fix a $p^{*} \in U_{\beta}$. For the sake of convenience, only in this proof, we will denote $A(S,\mathbf{q}_S,p^{*},\eta)$ by $A$, $B(S,p^{*},\eta)$ by $B$ and $O_{\beta}$ by $O$. Let $C$ be the event $\{X: \bar{p}_u(X) = \bar{q}_u \text{ for } u \in S\}$. Let $n$ be large enough so that the sets $A,B$ and $C$ are non-empty. Note that $\mathbbm{1}\{ X \in B\}$ is a function of $\Xmean$ and $\mathbbm{1}\{X \in C\}$ is a function of $\Xcav$. Therefore we write $\Xmean \in B$ and $\Xcav \in C$ in place of $X \in B$ and $X \in C$, respectively. 

With this notation in place, we have 
\begin{align}
    &\mu(A) = \frac{1}{Z_{\beta}}\sum_{\Xcav \in C}\sum_{\Xmean \in B}\exp(\mathcal{H}_{\beta}(X))\nonumber \\
    &\leq \frac{1}{Z_{\beta}}\sum_{\Xcav \in C}\sum_{\Xmean \in B}\exp(\Hmean(X;S))\exp\left(O(nr\eta + r^2) +  \Gamma(p^{*}) +\bar{\Delta}(\mathbf{q}_S,p^{*})\right) \nonumber \\
    &= \frac{H_{\mathsf{cav}}(\mathbf{q}_S)}{Z_{\beta}} \sum_{\Xmean \in B}\exp(\Hmean(X;S))\exp\left(O(nr\eta + r^2) +  \Gamma(p^{*}) +\bar{\Delta}(\mathbf{q}_S,p^{*})\right) \nonumber \\
    &\leq \frac{Z^{\mathsf{mean}}_{\beta}}{Z_{\beta}} \exp\Big(O(nr\eta + r^2) + 2n\sum_{u\in S}\Big[ L_{\beta}(\bar{q}_u)-L_{\beta}(p^{*})- \zeta_{\beta}(|\bar{q}_u-p^{*}|)\Big] \Big)\,,\label{eq:cavity_ub_1}
\end{align}
where 
$$\Gamma(p^{*}) := \sum_{i=0}^{K}\beta_i n r |V_i|(p^{*})^{|E_i|}\,,
 $$ 
 $$\bar{\Delta}(\mathbf{q}_S,p^{*}):= 2n\sum_{u\in S}\Big[ L_{\beta}(\bar{q}_u)+I(\bar{q}_u)-L_{\beta}(p^{*})-I(p^{*})- \zeta_{\beta}(|\bar{q}_u-p^{*}|)\Big] \,,$$ 
 and 
 $$Z_{\beta}^{\mathsf{mean}} := \sum_{\Xmean \in B} \exp(\Hmean(X;S))\exp(-2n|S|I(p^{*}) + \Gamma(p^{*}))\,.$$
 In the second step of \eqref{eq:cavity_ub_1} we have applied the cavity decomposition from Lemma~\ref{lem:cavity_decomp}. In the third step, we have used the fact that $|C| = H_{\mathsf{cav}}(\mathbf{q}_S)$ and in the fourth step, we have used Lemma~\ref{lem:count_approx} to upper bound this count.

 Now, note that under the condition $p^{*}< 1- |S|/n$, there exists an admissible restricted degree $\bar{q} \in \{0,1/n,\dots,1-{|S|}/{n}\}$ such that $|\bar{q} - p^{*}|\leq {1}/{n}$. Denote by $\mathbf{\bar{q}}$ the degree tuple with $\mathbf{\bar{q}}_u = \bar{q}$ for every $u \in S$.  Clearly, $\mu(B)\geq \mu(A(S,\mathbf{\bar{q}},p^{*},\eta))$. Repeating the calculation in Equation~\eqref{eq:cavity_ub_1}, but with corresponding lower bounds instead of upper bounds, we conclude that
  \begin{equation}\label{eq:mu_b_lb}
 \mu(B) \geq \frac{Z^{\mathsf{mean}}_\beta}{Z_\beta}\exp(-O(nr\eta+r^2\log(\tfrac{n}{r}) + r\log n) + \underline{\Delta}(\mathbf{\bar{q}},p^{*}))\,,
 \end{equation}
where
 $$\underline{\Delta}(\mathbf{q}_S,p^{*}) = 2n\sum_{u\in S}\Big[ L_{\beta}(\bar{q}_u)-L_{\beta}(p^{*})- C_{\beta}|\bar{q}_l-p^{*}|\Big]\,.$$

 Using the properties of the Shannon entropy (i.e, $\sup_{p \in [0,1-{1}/{n}]}|H(p)-H(p+{1}/{n})| = |H(0)-H({1}/{n})| \leq n\inv({1+ \log n})$), we have that %\gb{what property? easy bound on derivative?}, we have that 
 $$
 \sup_{p\in [0,1-1/n]}\Big|L_{\beta}(p)-L_{\beta}\Big(p+\frac{1}{n}\Big)\Big| \leq C_{\beta}\frac{\log n}{n}
 $$ 
 for some positive constant $C_{\beta}$. This implies that $\underline{\Delta}(\mathbf{\bar{q}},p^{*}) \leq C_{\beta}r\log(n)$. Plugging this into the lower bound on $\mu(B)$ in Equation~\eqref{eq:mu_b_lb}, and then combining it with the upper bound on $\mu(A)$ in Equation~\eqref{eq:cavity_ub_1}, we obtain the claim.
 \end{proof}

 From this result we derive the following corollary, which establishes that $\sup_{u\in [n]}|p_u(X)-p^{*}|$ must be close to zero with high probability under the measure $\mu\b(\cdot\bigr|\Ball \ps\eta\b)$. Observe that in the statement of Theorem~\ref{thm:cavity_conc}, the term in the exponent, $L_{\beta}(\bar{q}_u) - L_{\beta}(p^{*}) - \zeta_{\beta}(|\bar{q}_u-p^{*}|) < -\delta_0 < 0$ whenever $|\bar{q}_u-p^{*}|$ is large. Therefore, this the event where $|\bar{q}_u-p^{*}|$ is large incurs an exponentially small probability. We refer to Section~\ref{sec:cor_proof_1} for its complete proof. 
 
 %\gb{say a few words about why this follows.. accumulating probability... etc.}

 \begin{corollary}\label{cor:degree_conc2}\label{cor:degree_conc}
 Suppose $p^{*} \in U_{\beta}$ and let
$\epsilon > 0$ be an arbitrary fixed constant. Then, we can take $0<\eta < c(\beta,\epsilon) $ and $ n > n_0(\beta,\epsilon,\eta)$ such that
$$\mu\Big(\sup_{u \in [n]}|p_u(X) - p^{*}| \leq \epsilon\,\Big| \Ball \ps\eta\Big) \geq 1-\exp(-\Omega_{\beta,\epsilon}(n))\,.$$
 \end{corollary}

Recall that in Lemma~\ref{lem:almost_degree_conc}, we showed that whenever $X$ is close to the constant graphon $p^{*}$, then \emph{most} of the degrees concentrate close to $p^{*}$. In the result below we show that when $p_{u}(X)$ is close to $p^{*}$ uniformly for \emph{every} vertex $u$ (as shown in Corollary~\ref{cor:degree_conc}) and $X$ is close to the constant graphon $p^{*}$, most of the degrees $p_{uv}(X)$ concentrate close to $(p^{*})^2$. This will be useful in Section~\ref{sec:unif_conc}, where we will prove Theorem~\ref{thm:final_unif_conc}, which is the important component behind the proof of Theorem~\ref{thm:main_result}.

\begin{corollary}\label{cor:almost_dub_conc2}\label{cor:almost_dub_conc}
Suppose that $p^{*}\in U_{\beta}$. Given arbitrary $\epsilon, \delta > 0$, suppose $\eta < c(\beta,\epsilon,\delta)$ and $n > n_0(\delta,\epsilon,\beta,\eta)$. Then, for every $ u\in [n]$, there exists a random set $S_u \subseteq [n]\setminus \{u\}$ such that $|S_u|\leq \delta n$ and
$$
\mu\Big(\sup_{u\in [n]}\sup_{v\in S^{\complement}_u}|p_{uv}(X)-(p^{*})^2| \leq \epsilon\,\Big|\Ball \ps\eta\Big) \geq 1-\exp(-\Omega_{\beta,\epsilon}(n))\,.
$$
\end{corollary}
We refer to Appendix~\ref{sec:cor_proof_2} for the proof.

%% file: metaStabPf.tex
\section{Proof of Theorem~\ref{thm:local_metastability}}

\label{sec:ThmMetastable}

Before proceeding with the proof of Theorem~\ref{thm:local_metastability}, we will establish generalizations of \cite[Lemma 12 and Lemma 17]{bhamidi2011mixing}. Therefore, we will not instantiate to the model parameters given in Example~\ref{ex:t_ergm} but consider a general ERGM with parameter $\beta$. 
We need to treat the vertex $1$ separately from the other vertices. Following the notation preceding Lemma 12 in \cite{bhamidi2011mixing}, we define for some fixed, finite set of graphs $\mathcal{G}$:
$$\bar{r}_{\max}(X) := \max\Big(\max_{u:u\neq 1}p_u(X),\sup_{\substack{e = (u,v): u,v \neq 1\\ G \in \mathcal{G}}}r_G(X,e) \Big)$$
and
$$\bar{r}_{\min}(X) := \min\Big(\min_{u:u\neq 1}p_u(X),\inf_{\substack{e = (u,v): u,v\neq 1\\ G \in \mathcal{G}}}r_G(X,e)\Big) \,. $$
Here, we consider the evolution of the vertices $2,\dots,n-1$ when they are close to $G(n-1,p^{*})$ in terms of the subgraph counts and the edges connecting vertex $1$ are arbitrary. Notice that we have included the degrees $p_u(X)$ here in addition to $r_G(X,e)$, which will be useful to us later in the proof. The lemma below follows from a rewriting of the proof of Lemma 17 in \cite{bhamidi2011mixing}, by noting that the edges connected to vertex $1$ do not influence the evolution of $\delN e G X$ in the leading order term as considered in \cite[Lemma 12 and Lemma 14]{bhamidi2011mixing} and a straightforward tweak to also consider $p_u(X)$. Therefore, we skip the proof.  
\begin{lemma}\label{lem:non_exit_cavity}
Suppose $\mathcal{G} = \mathbb{G}_L$ (the set of all graphs with at-most $L$ vertices) and let $\epsilon > 0$ be a small enough constant independent of $n$. Suppose $p^{*} \in U_{\beta}$ and let $X_0,X_1,\dots$ are drawn from the Glauber dynamics with respect to the measure $\mu$. For some large enough $L \in \mathbb{N}$, independent of $n$, if $X_0$ is such that $p^{*}-\epsilon \leq \bar{r}_{\min}(X_0) \leq \bar{r}_{\max}(X_0) \leq p^{*} + \epsilon$, then for some $\alpha$ depending only on $\beta,\epsilon,L$, we have
$$\mathbb{P}\B(\sup_{t \leq e^{\alpha n}} \bar{r}_{\max}(X_t) \geq  p^{*} + 2\epsilon\B) \leq \exp(-\Omega_{\beta,L,\epsilon}(n))$$
and
$$\mathbb{P}\B(\inf_{t \leq e^{\alpha n}} \bar{r}_{\min}(X_t) \leq  p^{*} - 2\epsilon\B) \leq \exp(-\Omega_{\beta,L,\epsilon}(n)) \,.$$
\end{lemma}

We now instantiate our discussion to the case of the exponential random graph model defined in Example~\ref{ex:t_ergm} and use the notation established in this example. Recall $p_1(X)$ and $p_{1u}(X)$. Define $$
p^{(1)}_{\max}(X) = \max\Big(p_1(X),\max_{\substack{u\neq 1\\u\in[n]}}\frac{p_{1u}(X)}{p_1^{*}}\Big)\qquad\text{and}\qquad p^{(1)}_{\min}(X) = \min\Big(p_1(X),\min_{\substack{u\neq 1\\u\in[n]}}\frac{p_{1u}(X)}{p_1^{*}}\Big)\,.$$

\begin{lemma}\label{lem:metastable_cavity}
Consider the same setting as Lemma~\ref{lem:non_exit_cavity} instantiated to the parameter $\beta$ given in Example~\ref{ex:t_ergm}, with $p^{*} = p_1^{*}$. Given $\epsilon_1 > 0$, we can pick $\epsilon$ in Lemma~\ref{lem:non_exit_cavity} small enough such that the following holds:
Suppose $q^{*} -\epsilon_1 \leq p_{\min}^{(1)}(X_0)\leq p_{\max}^{(1)}(X_0) \leq q^{*} + \epsilon_1$,
then for some $\alpha$ depending only on $\beta,\epsilon,\epsilon_1,L$, we have
$$\mathbb{P}\B(\sup_{t \leq e^{\alpha n}} p_{\max}^{(1)}(X_t) \geq  q^{*} + 2\epsilon_1\B) \leq \exp(-\Omega_{\beta,L,\epsilon}(n)) $$
and
$$\mathbb{P}\B(\inf_{t \leq e^{\alpha n}} p^{(1)}_{\min}(X_t) \leq  q^{*} - 2\epsilon_1\B) \leq \exp(-\Omega_{\beta,L,\epsilon}(n)) \,.$$
\end{lemma}

\begin{proof}
Let $N := {{n}\choose{2}}$. Recall the function $g$ defined in Example~\ref{ex:t_ergm} and $\phi_{\beta}$ as defined in Section~\ref{sec:definition}. It is easy to show using similar techniques as in \cite[Lemma 12]{bhamidi2011mixing} that
$$-\frac{p_1(X_t)}{N} + \frac{g(p^{(1)}_{\min})}{N}\leq \mathbb{E}\left[p_1(X_{t+1})-p_1(X_t)|X_t\right] \leq -\frac{p_1(X_t)}{N} + \frac{g(p^{(1)}_{\max})}{N}
\,.$$

Similarly, for every $u \in [n]$ and $u \neq 1$, denoting $\bar{r}_{\max}(X_t), \bar{r}_{\min}(X_t)$ by $\bar{r}_{\max},\bar{r}_{\min}$ respectively,
$$ \mathbb{E}\left[p_{1u}(X_{t+1})-p_{1u}(X_t)|X_t\right] \leq -\frac{2p_{1u}(X_t)}{N} + \frac{g(p^{(1)}_{\max})\bar{r}_{\max}}{N} + \frac{p^{(1)}_{\max}\phi(\bar{r}_{\max})}{N}$$
and
$$\mathbb{E}\left[p_{1u}(X_{t+1})-p_{1u}(X_t)|X_t\right]  \geq -\frac{2p_{1u}(X_t)}{N} + \frac{g(p^{(1)}_{\min})\bar{r}_{\min}}{N} + \frac{p^{(1)}_{\min}\phi(\bar{r}_{\min})}{N}\,.$$
Now, notice that by Lemma~\ref{lem:non_exit_cavity}, $\bar{r}_{\max}(X_t)/p_1^{*} \leq 1 + 2\frac{\epsilon}{p_1^{*}}$ and $\bar{r}_{\min}(X_t)/p_1^{*} \geq 1 - 2\frac{\epsilon}{p_1^{*}}$ with probability at-least $\exp(-\Omega(n))$ whenever $t \leq \exp(\alpha n)$. Therefore, we can consider the evolution of $\frac{p_{1u}(X_t)}{p_1^{*}}$ akin to the evolution of $p_1(X_t)$ with $g(x)$ replaced by $\frac{g(x) + x}{2}$. Notice that the functions $g$ and $\frac{1}{2}(g(x)+x)$ play the role of $\phi()$ in the proof of \cite[Lemma 12]{bhamidi2011mixing} and satisfy the relationship $g(q^{*}) = q^{*}$ and $g^{\prime}(q^{*}) < 1$ (same for $\frac{g(x) + x}{2}$). This allows us to conclude the statement of the lemma with minor modifications to the proof of Lemma 17 in \cite{bhamidi2011mixing}
\end{proof}

\begin{proof}[Proof of Theorem~\ref{thm:local_metastability}]
Let the initial state $X_0$ be sampled as in the theorem statement.

\begin{enumerate}
    \item Recall $\delta_{\square}^{S,p^{*}}$ given in Definition~\ref{def:res_met}. Consider this with $S = \{1\}$. By Lemma~\ref{lem:restricted_metric}, in order to show that $X_t\in \Ball\ps\eta$ it is sufficient to show that $\delta_{\square}^{{1},p^{*}}(X_t,p^{*}) \leq {\eta}/{4}$ with high probability.  With similar arguments as in the proof of Lemma~\ref{lem:trajectory_inside_ball} with $r_{\min}, r_{\max}$ replaced with $\bar{r}_{\min},\bar{r}_{\max}$ we conclude that with probability at-least $1-T\exp(-\Omega_{\eta,\beta}(n))$, we have that every point in the trajectory $X_0,X_1,\dots,X_T \in \Ball\ps\eta$. Using Lemma~\ref{lem:metastable_cavity} and the result above we conclude the statement.
\item This follows from a straightforward application of Theorem~\ref{thm:cavity_conc} along with Theorem~\ref{thm:large_deviations} and the fact that $p_1^{*}$ is the unique global maximizer of $U_{\beta}$.

\item This follows from the same considerations as the proof of Item 2. 
\end{enumerate}
\end{proof}

%% file: techLems.tex
\section{Proofs of Technical Lemmas}
\label{sec:deferred_proofs}
\subsection{Proof of Lemma~\ref{lem:almost_degree_conc}}
% \begin{proof}
Fix any $p \in [0,1]$ and $\delta > 0$ and consider the sets $S^{+} := \bigr\{u: p_u(X)-p > {2\delta_{\square}(\tilde{X},p\mathbf{1}) }/{\delta}\bigr\}$ and $S^{-} := \{u: p_u(X)-p < -{2\delta_{\square}(\tilde{X},p\mathbf{1}) }/{\delta}\}$. By the definition of the cut-metric it follows that
$$\delta_{\square}(\tilde{X},p\mathbf{1}) \geq \frac{1}{n}\sum_{u\in S^{+}}p_u(X)-p \geq \frac{2|S^{+}|\delta_{\square}(\tilde{X},p\mathbf{1})}{n\delta}\,.
$$
Thus, $|S^{+}| \leq {\delta n}/{2}$ and similarly $|S^{-}| \leq {\delta n}/{2}$. Therefore, $|S^{+}\cup S^{-}| \leq \delta n$, which allows us to conclude the first inequality by taking $S = S^{+}\cup S^{-}$. The second inequality follows by directly applying Theorem~\ref{thm:large_deviations} to Equation~\eqref{eq:edge_cut_metric}.
% \end{proof}
\qed

\subsection{Proof of Lemma~\ref{lem:dom_bound_delta}}

We note that when considering distance from any constant graphon $p\mathbf{1}$ we have
$\delta_{\square}(\tilde{X},p\mathbf{1}) = \delta_{\square}(f^{X},p\mathbf{1})$, since measure-preserving operators do not affect the constant graphon. 
Now, notice that since $f^{Y}(x_1,x_2)-p\leq f^{X}(x_1,x_2)-p \leq f^{Z}(x_1,x_2)-p$, for any Borel measurable sets $S,T \subseteq [0,1]$ we have
$$\int_{S\times T} (f^{X}(x_1,x_2)-p)dx_1 dx_2 \leq \int_{S\times T} (f^{Z}(x_1,x_2)-p)dx_1 dx_2 \leq \delta_{\square}(\tilde{Z},p)\,. $$
Similarly,
$$\int_{S\times T} (f^{X}(x_1,x_2)-p)dx_1 dx_2 \geq \int_{S\times T} (f^{Y}(x_1,x_2)-p)dx_1 dx_2 \geq -\delta_{\square}(\tilde{Y},p)\,.$$
These together establish that
$$\B|\int_{S\times T} (f^{X}(x_1,x_2)-p)dx_1 dx_2\B| \leq \max(\delta_{\square}(\tilde{Y},p),\delta_{\square}(\tilde{Z},p))\,.$$
Taking the suprememum over $S,T$ yields the lemma. 
\qed

\subsection{Proof of Lemma~\ref{lem:young_ineq}}

In case $d_i < |E|$ for some $i$, we start with Young's product inequality, which states that for every $a,b \geq 0$ and $\alpha,\beta > 0$ such that ${\alpha}\inv+{\beta}\inv = 1$, we have $ab \leq {\alpha}\inv {a^{\alpha}}+{\beta}\inv{b^{\beta}}$. We take $a = q^{d_i}$, $b = p^{|E|-d_i}$, $\alpha = {|E|}/{d_i}$ and $\beta = {|E|}/{(|E|-d_i)}$. In case $d_i = |E|$, we check that the inequalities below hold trivially. Hence
\begin{align}
    \sum_{i=1}^{|V|}q^{d_i}p^{|E|-d_i} &\leq \sum_{i=1}^{|V|}\frac{d_i}{|E|}q^{|E|} + \frac{|E|-d_i}{|E|}p^{|E|}\nonumber \\
    &= \frac{\sum_{i=1}^{|V|}d_i}{|E|}q^{|E|} + \frac{|V||E|-\sum_{i=1}^{|V|}d_i}{|E|}p^{|E|}\nonumber \\
    &= 2 q^{|E|} + (|V|-2)p^{|E|}\,.
\end{align}
In the last step, we have used the fact that for any finite simple graph $G$, $\sum_{i=1}^{|V|}d_i = 2|E|$. Equality when $G = G_0$ follows by a straightforward calculation. 

Now suppose that $G \neq G_0$. Then, it is easy to show that there exists a vertex $j$ such that $d_j < |E|$. We note that Young's product inequality is strict whenever $a^{\alpha} \neq b^{\beta}$. For the choice of $a,b,\alpha,\beta$ above, this condition means $p \neq q$. Now, consider the function:
$f(p,q) =  -q^{d_j}p^{|E|-d_j} + \frac{d_i}{|E|}q^{|E|} + \frac{|E|-d_i}{|E|}p^{|E|}$. This is continuous over the set $[0,1]^2$. Define $A_{\delta} := \{(p,q) \in [0,1]^2 : |p-q|\geq \delta\}$. Clearly, $A_{\delta}$ is a compact set for every $\delta \geq 0$. Define $\zeta(\delta) = \inf_{(p,q) \in A_{\delta}} f(p,q) $. It is clear from the strictness of the Young's inequality that $f(p,q) > 0$ for every $(p,q) \in A_{\delta}$ whenever $\delta > 0$. Therefore, we conclude by compactness of $A_{\delta}$ and continuity of $f$ that $\zeta(\delta) > 0$ whenever $\delta > 0$. The continuity of $\zeta$ follows from the continuity of $f$. Therefore, we conclude that there exists $\zeta$ as in the statement of the lemma such that:
$$ -C|p-q|\leq q^{d_j}p^{|E|-d_j} - \frac{d_i}{|E|}q^{|E|} - \frac{|E|-d_i}{|E|}p^{|E|} \leq -\zeta(|p-q|)\,.$$ 
The inequality above holds with $\zeta = 0$ for every $i$, even when $d_i = |E|$. This allows us to sum the inequality above and conclude the result. 
\qed

\subsection{Proof of Lemma~\ref{lem:inner_count}}

Since the proof is elementary, we only provide a brief sketch. The statement $N_G^{0}(X,u) \geq N_G(X,u)$ follows from the fact that $A_u^{k} = \cup_{l=1}^{k}A^{k,l}_u$ and the union bound. Now, note that $0\leq \sum_{l=1}^{k} \mathbbm{1}(A^{k,l}_u) - \mathbbm{1}(A_u^{k}) \leq k$ and the sum is non zero only when the event $A_u^{k,l}\cap A_u^{k,m}$ holds for some $l\neq m$, $l,m \in [k]$. Noting that under the uniform measure over $[0,1]^k$ the measure of $A_u^{k,l}\cap A_u^{k,m}$ is ${n^{-2}}$ and using the union bound, we conclude the result. 
\qed

\subsection{Proof of Lemma~\ref{lem:single_vertex_decomp}}
% \begin{proof}
 Suppose $d_i$ is the degree of vertex $i \in [k]$. In light of Lemma~\ref{lem:inner_count} we may replace $N_G(X,u)$ in the lemma statement by $N_G^0(X;u)$ and then by considering a specific term in the sums we see that it is sufficient to prove that $$\biggr|\int_{[0,1]^{k}} \mathbbm{1}(A^{k,l}_u)\prod_{(i,j)\in E(G)}f^{X}(x_i,x_j) \prod_{i=1}^{k}dx_i - \frac{p_u(X)^{d_l}p^{|E(G)|-d_l}}{n}\biggr| \leq |E(G)|\frac{\delta_\square(\tilde{X},p)}{n}\,.$$
Notice that
\begin{align}
    &\int_{[0,1]^{k}} \mathbbm{1}(A^{k,l}_u)\prod_{(i,j)\in E(G)}f^{X}(x_i,x_j) \prod_{i=1}^{k}dx_i \nonumber \\ &= \int_{[0,1]^k}\mathbbm{1}(A^{k,l}_u)\prod_{(l,i)\in E(G)}f^{X}(x_l,x_i)\prod_{\substack{(i,j)\in E(G)\\ i,j\neq l}} f^{X}(x_i,x_j) \prod_{i=1}^{k}dx_i \,.
\end{align}

A simple computation shows that
\begin{align}
    \frac{p_u(X)^{d_l}p^{|E(G)|-d_l}}{n} &= \int_{[0,1]^k}\mathbbm{1}(A^{k,l}_u)\bigg(\prod_{\substack{i:(l,i)\in E(G)}}f^{X}(x_l,x_i) \bigg) p^{E(G)-d_l} \prod_{i=1}^{k}dx_i \,.
\end{align}
Therefore,
\begin{align}
&\biggr|\int_{[0,1]^{k}} \mathbbm{1}(A^{k,l}_u)\prod_{(i,j)\in E(G)}f^{X}(x_i,x_j) \prod_{i=1}^{k}dx_i - \frac{p_u(X)^{d_l}p^{|E(G)|-d_l}}{n}\biggr| \nonumber \\
 &=\biggr|\int_{[0,1]^{k}} \mathbbm{1}(A^{k,l}_u)\prod_{(l,i)\in E(G)}f^{X}(x_l,x_i)\bigg[\prod_{\substack{(i,j)\in E(G)\\ i,j\neq l}} f^{X}(x_i,x_j) - p^{|E(G)|-d_l}\bigg] \prod_{i=1}^{k}dx_i \biggr|
 \,.\label{eq:int_decomp_1}
\end{align}

Following the proof of \cite[Lemma 4.4]{borgs2006graph} with minor modifications, consider any ordering among the set of edges $(i,j) \in E(G)$ such that $i,j \neq l$ and index these ordered edges by $(i_1,j_1),\dots,(i_{h},j_{h})$ where $h = |E(G)| - d_l$.
Then
\begin{align*}&\bigg[\prod_{\substack{(i,j)\in E(G)\\ i,j\neq l}} f^{X}(x_i,x_j) - p^{|E(G)|-d_l} \bigg] \nonumber \\
&\qquad= \sum_{r = 0}^{|E(G)|-d_l-1}p^{r}\prod_{ m = r+1}^{|E(G)|-d_l}f^{X}(x_{i_m},x_{j_m}) - p^{r+1}\prod_{ m = r+2}^{|E(G)|-d_l}f^{X}(x_{i_m},x_{j_m})\,.  \end{align*}
Now, we use the above decomposition in Equation~\eqref{eq:int_decomp_1} and consider the terms in the summation one by one. We then follow the technique used in the proof of \cite[Lemma 4.4]{borgs2006graph} along with the fact that $A^{k,l}_u
$ depends only on $x_u$ and the fact that the measure of the event $A^{k,l}_u$ under the uniform measure over $[0,1]^k$ is ${1}/{n}$ to conclude the result. 
\qed

\subsection{Proof of Lemma~\ref{lem:inner_event_S}}
% \begin{proof}
 Only in this proof, we will take the probability space to be $[0,1]^k$ equipped with the Borel sigma algebra and the uniform measure $P$. 

First, note that by the union bound, 
$\sum_{u\in S}N_G(X;u) \geq N_G(X;S)$. Now, almost surely $$\sum_{u\in S} \mathbbm{1}(A^{k}_u) \leq k\,,$$ 
since there are at most $k$ vertices in the graph $G$. We conclude that almost surely $$\sum_{u\in S} \mathbbm{1}(A^{k}_u) -\mathbbm{1}(A^k_S)\leq k-1\,.$$
 Now let $B^k_S$ be the event that $\sum_l \mathbbm{1}(A^k_u) \neq \mathbbm{1}(A^{k}_S)$. This can happen only when two events $A^{k,i}_u$ and $A^{k,j}_v$ hold simultaneously for some $i,j \in [k]$, $u,v \in S$, $i\neq j$ and $u\neq v$. Therefore, we have 
 $$B^k_S = \bigcup_{\substack{u,v \in S \\ u\neq v}}\bigcup_{\substack{i,j\in [k]\\ i\neq j}} A_u^{k,i}\cap A_v^{k,j}\,.$$
 By the union bound,
$$P(B^{k}_S) \leq |S|^2k^2 P(A_l^{k,i}\cap A_m^{k,j}) = \frac{|S|^2k^2}{n^2}\,.$$
Now combining the considerations above, we have
\begin{align}
    &\sum_{l\in S}N_G(X;l) - N_G(X;S) \nonumber \\ &= \int_{[0,1]^k}\left[-\mathbbm{1}\left(A_S^{k}\right)+ \sum_{l\in S}\mathbbm{1}(A_l^{k})\right] \prod_{(i,j)\in E(G)} f^{X}(x_i,x_j)\prod_{i=1}^{k}dx_i \nonumber \\
    &\leq (k-1)\int_{[0,1]^k}\mathbbm{1}(B^k_S)  \prod_{(i,j)\in E(G)} f^{X}(x_i,x_j)\prod_{i=1}^{k}dx_i \nonumber \\
    &\leq k P(B^{k}_S) \leq \frac{|S|^2k^3}{n^2}\,.
\end{align} 
The lemma statement follows.
\qed
% \end{proof}
\subsection{Proof of Lemma~\ref{lem:restricted_metric}}
% \begin{proof}
In this proof alone, we will abuse notation to denote the set $\cup_{u \in S}[\frac{u-1}{n},\frac{u}{n}) \subseteq [0,1]$ also by $S$ (and similarly for $S^{\complement}$). Since we are considering the cut-metric between  $\tilde{X}$ and a constant graphon, we can write
 \begin{align}
     \delta_{\square}^{S,p^{*}}(\tilde{X},p^{*}) &= \sup_{A,B \subset [0,1]}\biggr|\int_{A\times B} \left[f^{X,S,p^{*}}(x_1,x_2)-p^{*}\right]dx_1 dx_2\biggr|
     \nonumber \\
     &= \sup_{A,B \subset [0,1]}\biggr|\int_{A\cap S^{\complement}\times B\cap S^{\complement}} \left[f^{X}(x_1,x_2)-p^{*}\right]dx_1 dx_2\biggr|
     \nonumber \\
     &= \sup_{A,B \subset S^{\complement}}\biggr|\int_{A\times B} \left[f^{X}(x_1,x_2)-p^{*}\right]dx_1 dx_2\biggr|\,.
 \end{align}
 From the last equality above, we conclude that $\delta_{\square}^{S,p^{*}}(\tilde{X},p^{*}) \leq \delta_{\square}(\tilde{X},p^{*})$. To conclude the lower bound, note that $$\delta_{\square}(\tilde{X},p) \leq \delta_{\square}(\tilde{X},\tilde{X}^{S,p^{*}}) + \delta_{\square}(\tilde{X}^{S,p^{*}},p^{*}) = \delta_{\square}(\tilde{X},\tilde{X}^{S,p^{*}}) + \delta^{S,p^{*}}_{\square}(\tilde{X},p^{*})\,.$$
It is now easy to show that $\delta_{\square}(\tilde{X},\tilde{X}^{S,p^{*}}) \leq 1 - {|S^{\complement}|^2}/{n^2} $, which when combined with the last display above proves the lower bound. 
% 
% \end{proof}\\
\qed

\subsection{Proof of Lemma~\ref{lem:count_approx}}
% \begin{proof}
By Stirling's formula, we have that for any $k \in \mathbb{N}$ and $p\in [0,1]$ such that $kp$ is an integer, 
$$\frac{\exp(-2kI(p))}{\sqrt{2k}}\leq {{k}\choose{kp}} \leq \exp(-2kI(p))\,.$$
A counting argument shows that $N_{\mathsf{cav}}(\mathbf{q}_{S}) = \prod_{u\in S}{{n-|S|}\choose{q_u n}}$. For the upper bound, note that
\begin{align}
    \prod_{u\in S}{{n-|S|}\choose{q_u n}} &\leq \prod_{u\in S}{{n}\choose{q_u n}} \leq \exp\B(-2n\sum_{u\in S}I(q_u)\B)\,.
\end{align}
Now, for the lower bounds, we note that whenever $q_u \leq q \leq 1$ and $q \geq 1/2$: $|2qI({q_u}/{q}) - 2I(q_u)| \leq (2+\log({1}/({1-q})))(1-q)$. Taking $q = ({n-|S|})/{n} = 1-{r}/{n}$ below, we have
\begin{align}
    \prod_{u\in S}{{n-|S|}\choose{q_u n}} &\geq  \frac{\exp\B(-\sum_{u\in S}2(n-|S|)I(\tfrac{nq_u}{n-|S|})\B)}{\left(\sqrt{2n}\right)^r}\nonumber \\
    &\geq \exp\B(-\sum_{u\in S}2nI(q_u) - r^2\left[4+2\log(\tfrac{n}{r})\right]-\frac{r}{2}\log(2n) \B)\,.
\end{align}
\qed
% \end{proof}

\section{Proofs of lemmas from Section~\ref{sec:unif_conc}}
\subsection{Proof of Lemma~\ref{lem:edge_remove_approx}}
Before proving the lemma, we derive an estimate for $ \delN e G X$.

Recall the event $A^{k}_u$ in \eqref{e:Aku} in Section~\ref{sec:cavity}. Suppose the fixed graph $G$ has the vertex set $[k]$. 
We now define for $u,v \in [n]$, $$A^{kG}_{uv} := \b\{x \in [0,1]^k: (\lfloor nx_i\rfloor,\lfloor nx_j\rfloor) \in \{(u,v), (v,u)\}  \text{ for some } (i,j) \in E(G) \b\}\,.$$
For $i,j \in [k]$ define
$$A^{kij}_{uv} := \b\{x \in [0,1]^k: (\lfloor nx_i\rfloor,\lfloor nx_j\rfloor) \in \{(u,v), (v,u)\} \b\}
\,.
$$

Now, the definition of homomorphism density yields
\begin{align}
    \delN e G X &= N_G(X^{+e}) - N_G(X^{-e}) \nonumber \\
    &= \int \prod_{(i,j)\in E(G)} f^{X^{+e}}(x_i,x_j)\prod_{t=1}^{k}dx_t - \prod_{(i,j)\in E(G)} f^{X^{-e}}(x_i,x_j)\prod_{t=1}^{k}dx_t \nonumber \\
    &= \int \mathbbm{1}(A^{kG}_{uv})\prod_{(i,j)\in E(G)} f^{X^{+e}}(x_i,x_j)\prod_{t=1}^{k}dx_t\,.
\label{eq:edge_add_1}
\end{align}
A computation similar to the proof of Lemma~\ref{lem:inner_count} shows that
     \begin{equation*}\label{eq:inner_edge_count}
     \int\B|\mathbbm{1}(A^{kG}_{uv}) - \sum_{(i,j)\in E(G)}\mathbbm{1}(A^{kij}_{uv})\B|\prod_{t=1}^{k}dx_t \leq \frac{E(G)^3}{n^3}\,.
     \end{equation*}  
     Using this in Equation~\eqref{eq:edge_add_1}, we conclude that
\begin{align}
    \delN e G X &= \int \sum_{(i,j) \in E(G)} \mathbbm{1}(A^{kij}_{uv})\prod_{(a,b)\in E(G)} f^{X^{+e}}(x_a,x_b)\prod_{t=1}^{k}dx_t  \pm O\B(\frac{E(G)^3}{n^3}\B) \nonumber \\
    &= \sum_{(i,j)\in E(G)}\int \mathbbm{1}(A^{kij}_{uv})\prod_{\substack{(a,b)\in E(G)\\ (a,b)\neq (i,j) }} f^{X^{+e}}(x_a,x_b)\prod_{t=1}^{k}dx_t  \pm O\B(\frac{E(G)^3}{n^3}\B) \label{eq:edge_diff_inner}\,.
\end{align}

We are now ready for the proof. 
\begin{proof}
We will prove the result by replacing $p_{uv}(X),p_u(X)$, and $\delta_{\square}(X,p^{*})$ with $p_{uv}(X^{+e}),p_u(X^{+e})$, and $\delta_{\square}(X^{+e},p^{*})$, noting that $\delta_{\square}(X,X^{+e})\leq {2}/{n^2}$, $|p_{uv}(X) - p_{uv}(X^{+e})| \leq {1}/{n}$, and $|p_u(X)-p_u(X^{+e})| \leq {1}/{n}$.

We will use $E$ and $E(G)$ interchangeably in this proof. For $(i,j)\in E(G)$, first consider the quantity
 \begin{align}
     &\int\mathbbm{1}(A^{kij}_{uv})\prod_{\substack{(a,b)\in E\\ (a,b)\neq (i,j) }} f^{X^{+e}}(x_a,x_b)\prod_{t=1}^{k}dx_t \nonumber \\ &= \int\mathbbm{1}(A^{kij}_{uv})\prod_{\substack{(i,l)\in E\\ l\neq j}}f^{X^{+e}}(x_i,x_l)\prod_{\substack{(j,l)\in E\\ l\neq i}}f^{X^{+e}}(x_j,x_l)\prod_{\substack{(a,b)\in E\\ a,b\notin \{i,j\} }} f^{X^{+e}}(x_a,x_b)\prod_{t=1}^{k}dx_t \nonumber \\
      &= (p^{*})^{|E|-d_{i}-d_j+1}\int\mathbbm{1}(A^{kij}_{uv})\prod_{\substack{(i,l)\in E\\ l\neq j}}f^{X^{+e}}(x_i,x_l)\prod_{\substack{(j,l)\in E\\ l\neq i}}f^{X^{+e}}(x_j,x_l)\prod_{t=1}^{k}dx_t \nonumber \\ &\quad\pm \frac{2|E(G)|\delta_{\square}(\tilde{X^{+e}},p^{*})}{n^2}\label{eq:decomp_3}\,.
 \end{align}
  In the last step we have used a similar peeling argument as in Lemma~\ref{lem:single_vertex_decomp}. Recalling the sets $A_u^{ki}$ from Section~\ref{sec:cavity}, we note that whenever $i \neq j$, $\mathbbm{1}(A^{kij}_{uv}) = \mathbbm{1}(A^{ki}_u\cap A^{kj}_v) + \mathbbm{1}(A^{kj}_u\cap A^{ki}_v)$. Now,
\begin{align}
    &\int\mathbbm{1}(A^{ki}_{u}\cap A^{kj}_v)\prod_{\substack{(i,l)\in E\\ l\neq j}}f^{X^{+e}}(x_i,x_l)\prod_{\substack{(j,l)\in E\\ l\neq i}}f^{X^{+e}}(x_j,x_l)\prod_{t=1}^{k}dx_t \nonumber \\
    &=\int\mathbbm{1}(A^{ki}_{u}\cap A^{kj}_v)\prod_{\substack{(i,l)\in E\\ l \neq j\\l \not\in  E_{ij}(G)}}f^{X^{+e}}(x_i,x_l)\prod_{\substack{(j,l)\in E\\ l \neq i \\l \in E_{ij}(G)}}f^{X^{+e}}(x_j,x_l) \times  \nonumber 
    \\ &\quad  
    \qquad \qquad\qquad\qquad\qquad\prod_{ l \in E_{ij}(G)}f^{X^{+e}}(x_i,x_l)f^{X^{+e}}(x_j,x_l)\prod_{t=1}^{k}dx_t \nonumber \\
    &= \frac{1}{n^2}(p_u(X^{+e}))^{d_i-1-d_{ij}}(p_v(X^{+e}))^{d_j-1-d_{ij}}(p_{uv}(X^{+e}))^{d_{ij}}\,.
\end{align}
In the last step, we have used the definitions of degrees $p_u$ and $p_{uv}$ given in Equations~\eqref{eq:norm_degree} and~\eqref{eq:norm_degree_2} in terms of integrals over $f^{X}$. Using this in Equation~\eqref{eq:decomp_3}, and shortening $p_u(X_e^{+}),p_v(X_e^{+})$ and $p_{uv}(X_e^{+})$ to $p_u,p_v$ and $p_{uv}$, we obtain
\begin{align}
    &n^2\int\mathbbm{1}(A^{kij}_{uv})\prod_{\substack{(a,b)\in E\\ (a,b)\neq (i,j) }} f^{X^{+e}}(x_a,x_b)\prod_{t=1}^{k}dx_t \nonumber\\&= (p^{*})^{|E|-d_i-d_j+1}(p_u)^{d_i-1-d_{ij}}(p_v)^{d_j-1-d_ij}(p_{uv})^{d_{ij}}\nonumber\\&\quad+(p^{*})^{|E|-d_i-d_j+1}(p_v)^{d_i-1-d_{ij}}(p_u)^{d_j-1-d_{ij}}(p_{uv})^{d_{ij}} \pm 2|E(G)|\delta_{\square}(\tilde{X^{+e}},p^{*})\,.
\end{align}
Using the condition $\sup_u|p_u(X) -p^{*}| < \epsilon $,  Equation~\eqref{eq:edge_diff_inner}, and the equation above, we conclude the statement of the lemma.
\end{proof}

\subsection{Proof of Lemma~\ref{lem:fixed_pt_conc}}

We first state a simple modification of \cite[Theorem 1.5]{chatterjee2007stein} below, which follows by essentially rewriting its proof with minor changes.

\begin{lemma}[Modification of Theorem 1.5 in \cite{chatterjee2007stein}]\label{lem:stein_modify}
Under the same notation as \cite[Theorem 1.5]{chatterjee2007stein}, we suppose the same conditions as the original theorem hold, except condition (ii), where we replace the assumption $\Delta(X) \leq Bf(X) +C$, with 
\begin{enumerate}
    \item $\Delta(X) \leq \alpha \mathbbm{1}(A) + \gamma \mathbbm{1}(A^{\complement})$ for some event $A \in \sigma(X)$, $\alpha,\gamma \in \mathbb{R}^{+}$
    \item $|f(X)| \leq M$ almost surely.
\end{enumerate} Then, we have
$$\mathbb{P}(|f(X)| > t) \leq \B(1+\frac{\gamma\exp(\theta_0 M)\mathbb{P}(A^{\complement})}{\alpha}\B)\inf_{\theta \in [0,\theta_0]}\exp\B(\frac{\theta^2\alpha}{2} - \theta t\B)\,.$$
\end{lemma}

\begin{proof}
In the proof of \cite[Theorem 1.5]{chatterjee2007stein}, in the display below Equation (7), we have for $\theta \in [0,\theta_0]$
\begin{align*}
    |m^{\prime}(\theta)| &\leq |\theta|\mathbb{E}\exp^{\theta f(X)}\Delta(X) \nonumber \\
    &\leq |\theta|\alpha m(\theta) + |\theta|\gamma \exp(|\theta| M)\mathbb{P}(A^{\complement}) \nonumber \\
    &\leq |\theta|\alpha m(\theta) + |\theta|\gamma \exp(|\theta_0| M)\mathbb{P}(A^{\complement}) \,.
\end{align*}
In the second step we have used the hypothesis that $\Delta(X) \leq \alpha \mathbbm{1}(A) + \gamma \mathbbm{1}(A^{\complement})$ and $|f(X)| \leq M$. Therefore
\begin{equation*}
    \frac{d}{d\theta} \log\left(m(\theta) + \tfrac{\gamma\exp(\theta_0 M)\mathbb{P}(A^{\complement})}{\alpha}\right) \leq \theta \alpha\,.
\end{equation*}
The result then follows by an application of Gronwall's lemma and the Chernoff bound.
\end{proof}

We will consider Glauber dynamics with respect to the measure $\mu(\spacedot|\Ball \ps{\eta})$ (Definition~\ref{def:glauber}) in order to generate the exchangeable pairs required by \cite[Theorem 1.5]{chatterjee2007stein} (and Lemma~\ref{lem:stein_modify}), where the event $\Ball \ps \eta$ is as defined in Equation~\eqref{eq:B_def}.  In the notation of \cite[Theorem 1.5]{chatterjee2007stein}, we consider $$F_{uv}(X,X^{\prime}):= \sum_{\substack{w \in [n]\\ w \notin \{u,v\}}}X_{uw}X_{vw} - \sum_{\substack{w \in [n] \\ w \notin \{u,v\}}}X^{\prime}_{uw}X^{\prime}_{vw}\,.$$

With the help of Lemma~\ref{lem:stein_modify}, we will now prove Lemma~\ref{lem:fixed_pt_conc}.

\begin{proof}[Proof of Lemma~\ref{lem:fixed_pt_conc}] By Lemma~\ref{lem:outer_prob}, we conclude that whenever $p^{*} \in U_{\beta}$ and $\eta > 0$ is small enough,
  $$\mu\B(\Ball \ps {\eta/2})\B|{\Ball \ps \eta }^{\complement}\B) \leq C(\beta,\eta)\exp(-c(\beta,\eta)n^2)\,.$$ 

Let $X$ be drawn from the distribution $\mu(\spacedot|\Ball \ps{\eta})$ and let $X^{\prime}$ be obtained by taking a single step of Glauber dynamics with respect to the measure $\mu(\spacedot|\Ball \ps{\eta})$. Clearly, $(X,X^{\prime})$ form an exchangeable pair. Let $$f_{uv}(X) := \mathbb{E}\left[F_{uv}(X,X^{\prime})|X\right]\,.$$ 
Notice that $|f_{uv}(X)| \leq 1$ almost surely since $|F_{uv}(X,X^{\prime})| \leq 1$. For $n$ large enough as a function of $\eta$, whenever $X \in \Ball \ps {\eta/2}$, the Glauber dynamics over $\mu(\bigr|{\Ball \ps \eta })$ is exactly equal to the Glauber dynamics w.r.t. $\mu$. Therefore, a simple calculation yields the following for $X \in \Ball \ps {\eta/2}$:
\begin{equation}\label{eq:cond_exp}
    f_{uv}(X) = \frac{4p_{uv}(X)}{n-1} -{{{n}\choose{2}}\inv}\sum_{\substack{w\in[n]\\ w\notin \{u,v\}}}\phi_{uw}(X_{\sim uw})X_{vw}+\phi_{vw}(X_{\sim vw})X_{uw} \,.
\end{equation}

For any edge $e$, we consider $\phi_{e}(X_{\sim e}) = \frac{\exp(\sum_{i=1}^{K}n^2\beta_i\delN {e} i X )}{1+\exp(\sum_{i=1}^{K}n^2\beta_i\delN {e} i X )}$. Note that
% If $X_{\sim e} = X_{\sim e}^{\prime}$, then $|\phi_{e}(X_{\sim e}) - \phi_{e}(X^{\prime}_{\sim e})| = 0$. 
% 
\begin{equation}\label{eq:lips_probab}|\phi_{e}(X_{\sim e}) - \phi_{e}(X^{\prime}_{\sim e})| \leq C(\beta)n^2\sup_{i \in [K]}| \delN {e} i X -\delN {e} i {X^{\prime}}|\,.
\end{equation}
Now, $X_{\sim e}$ and $X_{\sim e}^{\prime}$ can differ at most in one edge, by construction. Suppose this edge is $h$. When $e = h$, then $\delN {e} i X  = \delN {e} i {X^{\prime}} $. Now suppose $e \neq h$. 
Invoking Equation~\eqref{eq:edge_add_1} with $k = |V_i|$, we obtain
\begin{align}
    &|\delN e i X - \delN e i {X^{\prime}}| \nonumber \\
    &= \biggr|\int \mathbbm{1}(A^{kG_i}_{e})\B[\prod_{(a,b)\in E(G)} f^{X_e^{+}}(x_a,x_b)-\prod_{(a,b)\in E(G)} f^{(X^{\prime}_e)^{+}}(x_a,x_b)\B]\prod_{t=1}^{k}dx_t \biggr| \nonumber \\
    &\leq \int \mathbbm{1}(A^{kG_i}_{e})\mathbbm{1}(A^{kG_i}_h)\prod_{t=1}^{k}dx_t \leq \frac{|E(G_i)|^2}{n^3}\,. \label{eq:change_count}
\end{align}
Combining Equations~\eqref{eq:cond_exp},~\eqref{eq:lips_probab} and~\eqref{eq:change_count}, we conclude that whenever $X \in \Ball \ps {{\eta}/{2}}$, \begin{equation}\label{eq:lip_func}
    |f_{uv}(X) - f_{uv}(X^{\prime})| \leq \frac{C(\beta)}{n^2}\,.
\end{equation}

Observe that whenever the Glauber dynamics does not update an edge of the form $(u,w)$ or $(v,w)$, $F_{uv}(X,X^{\prime}) = 0$. Let $A^{\mathsf{upd}}$ denote the event where $F_{uv}(X,X^{\prime}) \neq 0$. Clearly, $\mathbb{P}(A^{\mathsf{upd}}|X) \leq {4}/{n}$. It is also clear that $|F_{uv}(X,X^{\prime})| \leq 1$ almost surely. Therefore, we have for any $X \in \Omega$
\begin{equation}\label{eq:crude_bnd}
    |f(X)| \leq \mathbb{E}\left[|F(X,X^{\prime})||X\right] \leq \mathbb{P}\b(A^{\mathsf{upd}}|X\b) \leq \frac{4}{n}\,.
\end{equation}

Now consider the local variance proxy $\Delta_{uv}(X)$ (where the notation is once again derived from \cite[Theorem 1.5]{chatterjee2007stein}) whenever $X \in \Ball \ps{\eta/2}$,
\begin{align}
 \Delta_{uv}(X) &:= \frac{1}{2}\mathbb{E}\left[\left(f_{uv}(X)-f_{uv}(X^{\prime})\right)F_{uv}(X,X^{\prime})\bigr|X\right] \nonumber \\
 &\leq \frac{C(\beta)}{n^2}\mathbb{E}\b[\mathbbm{1}(A^{\mathsf{upd}})\bigr|X\b] \leq \frac{C(\beta)}{n^3}\,.
\end{align}
Here, we have used Equation~\eqref{eq:lip_func}. Whenever, $X \notin \Ball \ps{\eta/2}$, we will use the crude bound $$\Delta_{uv}(X) \leq \frac{C}{n^2}\,,$$
obtained by plugging in $|f(X)-f(X^{\prime})| \leq |f(X)| + |f(X^{\prime})| \leq {8}/{n}$ (which follows from Equation~\eqref{eq:crude_bnd}) into the the definition of $\Delta_{uv}$. Combining these bounds, we get that
$$\Delta_{uv}(X) \leq \frac{C(\beta)}{n^3}\mathbbm{1}(X \in \Ball \ps{\eta/2}) + \frac{C(\beta)}{n^2}\mathbbm{1}(X \not\in \Ball \ps{\eta/2})\,.$$

By Lemma~\ref{lem:outer_prob},
$$\mu\B(X \not\in \Ball \ps{\eta/2}\B|\Ball \ps{\eta/2} \B) \leq C_{\beta,\eta}\exp(-c(\beta,\eta)n^2)\,.$$
Now, applying Lemma~\ref{lem:stein_modify}, with $M = {4}/{n}$, and $\theta_0 = C_\beta n^2$ for some large enough $C_{\beta}$ and $t = {4\gamma}/({n-1})$, we conclude the result.
\end{proof}

\section{Deferred Proofs for Path Coupling}
\label{sec:deferred2}

%\gb{need to check that we're OK (i.e., adjust all proofs slightly) without max and min versions of $r_G$, or else define here}

\subsection{Proof of Lemma~\ref{lem:approx_mix}}
% \begin{proof}
Consider the following coupling between the trajectories $X_0,\dots,X_K$ and $Y_0,\dots,Y_K$:
\begin{enumerate}
    \item Generate $(X_0,Y_0)$ from the specified initial distribution.
    \item Given $X_k,Y_k$, we generate $X_{k+1},Y_{k+1}$ as follows:
    \begin{equation}
        (X_{k+1},Y_{k+1})\sim \begin{cases} Q_{X_k,Y_k} &\quad \text{ if } (X_k,Y_k) \in A \\
        P(X_k,\cdot)\times P(Y_k,\cdot) &\quad \text{ otherwise}
        \end{cases}
    \end{equation} 
\end{enumerate}
Now we consider the distance $d(X_{k+1},Y_{k+1})$. We have
\begin{align}
    \mathbb{E}d(X_{k+1},Y_{k+1}) &= \mathbb{E}d(X_{k+1},Y_{k+1})\mathbbm{1}((X_k,Y_k) \in A\times A)\nonumber \\ &\quad + \mathbb{E}d(X_{k+1},Y_{k+1})\mathbbm{1}((X_k,Y_k)\in A^{\complement}) \nonumber \\
     &\leq \mathbb{E}(1-\gamma)d(X_{k},Y_{k})\mathbbm{1}((X_k,Y_k) \in A) \nonumber \\ &\quad + \mathbb{E}d(X_{k+1},Y_{k+1})\mathbbm{1}((X_k,Y_k) \in A^{\complement}) \nonumber \\
    &\leq (1-\gamma)\mathbb{E}d(X_{k},Y_{k}) + \mathbb{E}d(X_{k+1},Y_{k+1})\mathbbm{1}((X_k,Y_k)\in A^{\complement}) \nonumber \\
    &\leq (1-\gamma)\mathbb{E}d(X_{k},Y_{k}) + \bar{D}\mathbbm{P}((X_k,Y_k)\in A^{\complement}) \nonumber \\ 
    &\leq (1-\gamma)\mathbb{E}d(X_{k},Y_{k}) + \bar{D}p_k\,.
\end{align}
We conclude the result by unrolling the recursion.
% \end{proof}
\qed

\subsection{Proof of Lemma~\ref{lem:stoc_dom}}

In both this and the next proof, we let $$r_{\min}(X) := \inf_{\substack{e\in {{n}\choose{2}}\\ G \in \mathbb{G}_L}} r_{G}(X,e)\qquad \text{and} \qquad r_{\max}(X) := \sup_{\substack{e\in {{n}\choose{2}}\\ G \in \mathbb{G}_L}} r_{G}(X,e)\,.$$

\begin{proof}
We will only prove the coupling for $G(n,p^{*}+\epsilon)$. The other coupling follows analogously. Let $N = {{n}\choose{2}}$. Let $e_1,\dots,e_N$ be any enumeration of ${{[n]}\choose{2}}$. Let $X_0 \sim \mu(\spacedot|\Ball \ps{\eta})$ and obtain the sequence $X_0,\dots,X_N$ by updating as follows: Given $X_{i-1}$, define $(X_{i})_e = (X_{i-1})_{e}$ whenever $e \neq e_i$ and let $(X_{i})_{e_i}$ be independently re-sampled from the conditional distribution $(X_{i-1})_{e_i}\bigr|(X_{i-1})_{\sim e_i}$. In other words, we obtain $X_{i}$ by re-sampling the coordinate $e_i$ in $X_{i-1}$. Clearly, $X_N \sim \mu(\spacedot|\Ball \ps{\eta})$.

Now, consider $Y_0 = 0 \in \Omega$ almost surely. Let $e_1,\dots,e_N$ be the same as before. We will construct $Y_1,\dots,Y_N$ as follows: Given $Y_{i-1}$, we construct $Y_{i}$ such that $(Y_i)_{\sim e_i} = (Y_{i-1})_{\sim e_i}$ and $(Y_i)_{e_i}$ is freshly drawn from $\mathsf{Ber}(p^{*} + \epsilon)$.  It is clear that $Y_N \sim G(n,p^{*}+\epsilon)$.

By Theorem~\ref{thm:final_unif_conc}, we have that $ p^{*} -\epsilon \leq r_{\min}(X_i)\leq r_{\max}(X_i) \leq p^{*}+\epsilon$ with probability at-least $1-C(\eta,\beta,\epsilon)\exp(-c(\eta,\beta,\epsilon)n)$. Recall the definition of $\phi_{\beta}$ from Section~\ref{sec:definition} and note that  $\mathbb{P}((X_{i})_{e_i}=1|X_{i-1}) \leq \phi_{\beta}(r_{\max}(X_{i-1})) \leq p^{*}+\epsilon$. Therefore the fresh draws can be coupled such that $(X_{i})_{e_i} \leq (Y_i)_{e_i}$ with probability at least $1-C(\eta,\beta,\epsilon)\exp(-c(\eta,\beta,\epsilon)n)$. Since $(X_i)_{e_i} = (X_N)_{e_i}$ and $(Y_i)_{e_i} = (Y_N)_{e_i}$, we conclude via a union bound over $i \leq N$ that $X_N \preceq Y_N$ with probability at least $1-C(\eta,\beta,\epsilon)\exp(-c(\eta,\beta,\epsilon)n)$. This gives the desired coupling by taking $X= X_N$ and $\bar{Y} = Y_N$ in the statement of the lemma.
\end{proof}

\subsection{Proof of Lemma~\ref{lem:trajectory_inside_ball}}
    % \begin{proof}
    We first show that with probability at-least $1-C(\eta)\exp(-c(\eta)n)$, we must have $p^{*} - {\eta}/{4}\leq r_{\min}(X_0)\leq r_{\max}(X_0) \leq p^{*} + {\eta}/{4}$.
    
    This can be shown for example by using Lemma~\ref{lem:edge_remove_approx} and simple concentration bounds for the degrees $p_u(X_0)$ and $p_{uv}(X_0)$. The fact that $\delta_{\square}(\tilde{X}_0,p^{*})$ is small follows from Theorem~\ref{thm:large_deviations} since $G(n,p^{*})$ is also (a very special case of) an exponential random graph. Now, invoking Lemma~\ref{lem:non_exit}, we conclude that with probability at least $1-TC(\beta,\eta)\exp(-c(\beta,\eta)n)$, we have for every $t \leq T$ that
    \begin{equation}\label{eq:r_bounded_prob}
    p^{*} - \frac{\eta}{2}\leq r_{\min}(X_t)\leq r_{\max}(X_t) \leq p^{*} + \frac{\eta}{2}\,.
    \end{equation}
    
    Now consider Markov chains $Y_0,\dots,Y_T$ and $Z_0,\dots,Z_T$ where $Y_0 \sim G(n,p^{*}-{\eta}/{2})$ and $Z_{0}\sim G(n,p^{*}+\frac{\eta}{2})$. Here, we generate the respective trajectories by Glauber dynamics with respect to $G(n,p^{*}-{\eta}/{2})$ (resp. $G(n,p^{*}+\frac{\eta}{2})$ ) We couple the trajectories as follows:
    \begin{enumerate}
        \item At step $0$, we pick $X_0,Y_0,Z_0$ such that almost surely
        $$ Y_0 \preceq X_0 \preceq Z_0\,.$$
        \item At step $t$, pick the same edge $E_t \sim \mathsf{unif}\b({{[n]}\choose{2}}\b)$ to update for each $X_{t-1},Y_{t-1}$ and $Z_{t-1}$.
        \item Pick $u_t \sim \mathsf{unif}([0,1])$ independently of everything else and set
        $$(X_{t})_{E_t} = \mathbbm{1}(u_t \leq \phi_{E_t}((X_t)_{\sim E_t})); \quad (Y_{t})_{E_t} = \mathbbm{1}(u_t \leq p^{*}-\tfrac{\eta}{2});\quad (Z_{t})_{E_t} = \mathbbm{1}(u_t \leq p^{*}+\tfrac{\eta}{2})\,.$$
    \end{enumerate}
    
    For $\eta$ small enough, we verify that under the event in Equation~\eqref{eq:r_bounded_prob} $$p^{*}-\frac{\eta}{2}\leq \phi_{\beta}(r_{\max}(X_t))\leq \phi_{E_t}((X_t)_{\sim E_t}) \leq \phi_{\beta}(r_{\max}(X_t)) \leq p^{*}+ \frac{\eta}{2}\,.$$
    This implies that $(Y_t)_{E_t} \leq (X_t)_{E_t} \leq (Z_t)_{E_t}$. We conclude that with probability at least $1-TC(\beta,\eta)\exp(-c(\beta,\eta)n)$ we have $
    Y_t \preceq X_t \preceq Z_t$. Now, we will apply Theorem~\ref{lem:dom_bound_delta} to obtain \begin{equation}\label{eq:graphon_bd_1}
    \delta_{\square}(\tilde{X}_t,p^{*})\leq \max\b(\delta_{\square}(\tilde{Y}_t,p^{*}),\delta_{\square}(\tilde{Z}_t,p^{*})\b)\,.
    \end{equation}
    Now, since $Y_t \sim G(n,p^{*}-{\eta}/{2})$ and $Z_t \sim G(n,p^{*}+{\eta}/{2})$, we have
    $$\mathbb{P}\B(\max\b(\delta_{\square}(\tilde{Y}_t,p^{*}),\delta_{\square}(\tilde{Z}_t,p^{*})\b) > \eta\B) \leq 1-C(\eta)\exp(-c(\eta)n^2) \,.$$
    Using this in Equation~\eqref{eq:graphon_bd_1}, we conclude the statement of the claim.
    % \end{proof}
    \qed
    
\section{Proof of Corollary~\ref{cor:degree_conc}} 
\label{sec:cor_proof_1}
\begin{proof}
Let $\delta > 0$ and $\eta >0$. We are given $X \in \Omega$ such that $X \in \Ball \ps\eta$, i.e., $\delta_{\square}(\tilde{X},p^{*}) < \eta$. Consider the set $S(X) := \{u: |p_u(X) - p^{*}| > {2\eta}/{\delta}\}$. By Lemma~\ref{lem:almost_degree_conc}, we note that $|S(X)| \leq \delta n$, so
 $$\{S(X)\neq \emptyset\} = \bigcup_{\substack{S \subseteq [n]\\1\leq |S|\leq \delta n}} \{S(X) = S\}\,.$$
By the union bound, 
 \begin{align}
     \mu\b(|S(X)|\neq \emptyset\bigr|\Ball \ps\eta\b)&\leq \sum_{r=1}^{\delta n}\sum_{\substack{S \subseteq [n]\\|S| = r}}\mu\b(S(X) = S|\Ball \ps\eta\b) \nonumber \\
     &= \sum_{r=1}^{\delta n}{{n}\choose{r}}\mu\b(S(X) = [r]|\Ball \ps\eta\b)\,.\label{eq:union_bd_1}
 \end{align}
 Here, the second step follows from the permuation invariance of the vertices with respect to the measure $\mu$. In order to evaluate the upper bound in Equation~\eqref{eq:union_bd_1}, we will consider the measure of the event $\{S(X) = [r]\}$ for $r \leq \delta n$. 
 
Consider the restricted degrees $\bar{p}_u$ with respect to the set $[r]$ and suppose that $\delta^2 < 2\eta$. Let the set $D(p^{*},\alpha,S) := \{\mathbf{q}_S: \inf_{u \in S}|\bar{q}_u - p^{*}| > \alpha\}$ and note that
 \begin{align}\label{eq:decomp_1}
     \{S(X) = [r]\} &= \{|p_u(X)-p^{*}|> {2\eta}/{\delta}, \forall u \in [r]\}  \nonumber \\
     &\subseteq \{|\bar{p}_u(X)-p^{*}|> {2\eta}/{\delta}-\delta, \forall u \in [r]\} \nonumber \\
     &= \bigcup_{\mathbf{q}_{[r]} \in D(p^{*},\frac{2\eta}{\delta}-\delta,[r])} \{\bar{p}_u(X) = \bar{q}_u, \forall u \in [r]\}\,.
 \end{align}
 Therefore whenever $\delta,\eta$ satisfy ${\epsilon}/{2} = {2\eta}/{\delta} -\delta$, we have
 \begin{align}
     \{S(X) = [r]\}\cap \{X \in B([r],p^{*},\eta)\} \subseteq \bigcup_{\mathbf{q}_{[r]} \in D(p^{*},\tfrac{\epsilon}{2},[r])}A([r],\mathbf{q}_{[r]},p^{*},\eta)
     \,.
     \label{eq:decomp_2}
 \end{align}
 
 Note that with the above choice of $\eta$ in the definition of $S(X)$, we can take $\delta < {\epsilon}/{2}$ to conclude that for all $u \in S(X)^{\complement}$, $|p_u(X)-p^{*}|\leq \epsilon$.
 We will note a simple result which follows from standard arguments in calculus.
 
 \begin{lemma}\label{lem:basic_calc}
 Suppose $p^{*} \in U_{\beta}$ and $q \in [0,1]$ is such that $|q-p^{*}| > \epsilon/2$. Then, $L_{\beta}(q)-L_{\beta}(p^{*}) - \zeta_{\beta}(|q-p^{*}|) < -C(\beta,\epsilon) < 0$.
 \end{lemma}
 
Pick $\delta$ to be small enough such that $p^{*} \in (\delta,1-\delta)$. Combining Lemma~\ref{lem:basic_calc} with Equation~\eqref{eq:decomp_2} and Theorem~\ref{thm:cavity_conc} we conclude that whenever $r \leq n\delta$:
 \begin{align}
     &\frac{\mu(\{S(X) = r\}\cap B([r],p^{*},\eta))}{\mu(B([r],p^{*},\eta))} \leq \sum_{\mathbf{q}_{[r]} \in D(p^{*},\tfrac{\epsilon}{2},[r])} \frac{\mu(A([r],\mathbf{q}_{[r]},p^{*},\eta))}{\mu(B([r],p^{*},\eta))} \nonumber \\
     &\leq |D(p^{*},\tfrac{\epsilon}{2},[r])|\exp(-nrC(\beta,\epsilon) + O_{\beta}(nr\eta + r^{2}\log(\tfrac{n}{r}) + r\log n)) \nonumber \\
     &\leq \exp(-2nrC(\beta,\epsilon) + O_{\beta}(nr\eta + r^{2}\log(\tfrac{n}{r}) + r\log n))\nonumber \\
     &= \exp\left(-nr\left(2C(\beta,\epsilon)- O_{\beta}(\eta + \delta\log(\tfrac{1}{\delta})+ \tfrac{\log n}{n})\right)\right)\nonumber \\
     &= \exp\left(-\Omega_{\beta,\epsilon}(nr) \right) \,.
 \end{align}
The first step follows from the union bound on Equation~\eqref{eq:decomp_2} and in the second step we used Theorem~\ref{thm:cavity_conc} along with Lemma~\ref{lem:basic_calc}. In the third step, we used that $|D(p^{*},\epsilon,[r])| \leq n^r$. In the last step, we have picked $\eta, \delta$ small enough (as functions of $\epsilon,\beta$) so that ${2\eta}/{\delta}-\delta = {\epsilon}/{2}$, $\delta < {\epsilon}/{2}$, $\eta < \eta_0(\epsilon,\beta)$ and $n$ large enough as a function of $\epsilon,\beta$ such that $O_{\beta}(\eta + \delta\log(1/{\delta})+ {n}\inv{\log n})) < C(\beta,\epsilon)$.

  By Lemma~\ref{lem:restricted_metric}, whenever $r\leq \delta n$, we must have $\Ball \ps\eta=B(\emptyset,p^{*},\eta)\subseteq B([r],p^{*},\eta) \subseteq B(\emptyset,p^{*},\eta+\delta)$. Along with Lemma~\ref{lem:outer_prob}, we conclude that whenever $\eta,\delta$ are smaller than some constant $c_{\beta} > 0$, we have
 $$\mu(\Ball \ps \eta) \geq c(\beta,\eta,\delta)\mu(B([r],p^{*},\eta))\,.$$ 
 Therefore whenever $n$ is larger than a constant depending only on $\beta,\epsilon$ and $\eta_0$, we have 
 \begin{align}
     \frac{\mu(\{S(X) = [r]\}\cap B(\emptyset,p^{*},\eta))}{\mu(B(\emptyset,p^{*},\eta))} &\leq \frac{\mu(\{S(X) = r\}\cap B([r],p^{*},\eta))}{\mu(B(\emptyset,p^{*},\eta))} \nonumber \\
     &\leq C(\beta,\eta,\delta)\frac{\mu(\{S(X) = r\}\cap B([r],p^{*},\eta))}{\mu(B([r],p^{*},\eta))}\nonumber \\
     &\leq \exp(-\Omega_{\beta,\epsilon}(nr))\,. \label{eq:right_cond}
 \end{align} 
 Using the fact that ${{n}\choose{r}} \leq n^r$, it follows using Equations~\eqref{eq:union_bd_1} and~\eqref{eq:right_cond} that
 \begin{align*}
     \mu\b(|S(X)|\neq \emptyset\bigr|\Ball \ps \eta\b) &\leq \sum_{r=1}^{n\delta} n^r \exp(-\Omega_{\beta,\epsilon}(nr))\\
     &\leq \exp(-\Omega_{\beta,\epsilon}(n))\,,
 \end{align*}
which proves the corollary.
\end{proof}

\subsection{Proof of Corollary~\ref{cor:almost_dub_conc}}
\label{sec:cor_proof_2}
% \begin{proof}

Recall the sets $A^{k,i}_u := \{x \in [0,1]^k: x_i \in  [\frac{u-1}{n},\frac{u}{n}) \}$. Let $\eta > 0$ be arbitrary for now. Suppose that $X \in \Ball \ps \eta$. We have by definition
$$p_{uv}(X) = n^2 \int f^X(x_1,x_3)f^X(x_2,x_3)\mathbbm{1}(A^{3,1}_u)\mathbbm{1}(A^{3,2}_v)dx_1 dx_2 dx_3\,.$$
A simple calculation reveals that
$$p_{uv}(X) - p_u(X)p^{*} = n^2\int f^X(x_1,x_3)\left(f^X(x_3,x_2)-p^{*}\right)\mathbbm{1}(A^{3,1}_u)\mathbbm{1}(A^{3,2}_v)dx_1 dx_2 dx_3\,.$$
Now, define $S_u^{+} := \{v: p_{uv}(X) - p^{*}p_u(X) > {\epsilon}/{2}\}$ and $S_u^{-} := \{v: p_{uv}(X) - p^{*}p_u(X) <-{\epsilon}/{2}\}$. Summing the display above for $v \in S^{+}_u$, we have
\begin{align}
    \frac{|S_u^{+}|\epsilon}{2n^2}  &< \sum_{v \in S_u^{+}}\int f^X(x_1,x_3)\left(f^X(x_3,x_2)-p^{*}\right)\mathbbm{1}(A^{3,1}_u)\mathbbm{1}(A^{3,2}_{v})dx_1 dx_2 dx_3 \nonumber \\
    &= \int_{x_2 \in A_u^{+}} f^X(x_1,x_3)\left(f^X(x_3,x_2)-p^{*}\right)\mathbbm{1}(A^{3,1}_u) dx_1 dx_2 dx_3 \nonumber \\
    &= \int_{x_1 \in [\tfrac{u-1}{n},\tfrac{u}{n})}\B[\int_{\substack{x_2 \in A_u^{+}\\ x_3 \in C(x_1)}} \left(f^X(x_3,x_2)-p^{*}\right)dx_2 dx_3\B] dx_1 \nonumber \\
    &\leq \frac{\delta_{\square}(\tilde{X},p^{*})}{n}\,,
\end{align}
where $A_u^{+} :=  \cup_{v\in S_u^{+}}[\tfrac{v-1}{n},\tfrac{v}{n})$ and  $C(x_1) := \{x \in [0,1]: f^{X}(x_1,x) = 1\}$. 

Taking $\eta < {\delta\epsilon}/{2}$, we conclude that
$$|S_u^{+}| \leq \frac{2n \delta_{\square}(\tilde{X},p^{*})}{\epsilon} \leq \frac{\delta n}{2} $$
and similarly
$$|S_u^{-}| \leq \frac{2n \delta_{\square}(\tilde{X},p^{*})}{\epsilon} \leq \frac{\delta n}{2}\,. $$
Therefore, whenever $\delta_{\square}(\tilde{X},p^{*}) < \eta < {\delta \epsilon}/{2}$, the sets $S_u := S_u^{+}\cup S_{u}^{-}$ are such that
\begin{equation}\label{eq:without_vert}
    \sup_{u\in [n]}\sup_{v\in S^{\complement}_u}|p_{uv}(X)-p^{*}p_u(X)| \leq \frac{\epsilon}{2}\,.
\end{equation}

Now, invoking Corollary~\ref{cor:degree_conc}, we conclude that for any $\eta < c(\beta,\epsilon,\delta)$ and whenever $n$ is larger than a constant depending only on $\beta, \epsilon,\eta$ and $\delta$ we have
$$\mu\B(\sup_{u \in [n]}|p_u(X) - p^{*}| \leq \frac{\epsilon}{2}\B| \Ball \ps \eta \B) \geq 1-\exp(-\Omega_{\beta,\epsilon}(n))$$
and from Equation~\eqref{eq:without_vert}, we have
$$\mu\B(\sup_{u\in [n]}\sup_{v\in S^{\complement}_u}|p_{uv}(X)-p^{*}p_u(X)| \leq \frac{\epsilon}{2}\B| \Ball \ps \eta\B) = 1\,.$$
Combining the two displays above, the statement of the corollary follows.
% \end{proof}
\qed